\documentclass{amsproc}
\usepackage{amssymb}

\newcommand{\Alt}[2]{Alt{\i}nok$_#1(#2)$}
\newcommand{\rd}[1]{\lfloor #1\rfloor}
\newcommand{\half}{\frac12}
\newcommand{\recip}[1]{\frac1{#1}}
\newcommand{\Span}[1]{\left<#1\right>}
\newcommand{\fie}{\varphi}
\newcommand{\la}{\lambda}
\newcommand{\si}{\sigma}
\newcommand{\ep}{\varepsilon}
\newcommand{\Ga}{\Gamma}
\newcommand{\La}{\Lambda}
\newcommand{\Si}{\Sigma}
\newcommand{\bmu}{\mathbf{\mu}}
\newcommand{\LK}{\La_{\mathrm{K3}}} % the K3 lattice
\newcommand{\C}{\mathbb C}
\newcommand{\FF}{\mathbb F}
\newcommand{\PP}{\mathbb P}
\newcommand{\Q}{\mathbb Q}
\newcommand{\Z}{\mathbb Z}
\newcommand{\Oh}{\mathcal O}
\newcommand{\sB}{\mathcal B}
\newcommand{\sC}{\mathcal C}
\newcommand{\sE}{\mathcal E}
\newcommand{\sL}{\mathcal L}
\newcommand{\broken}{\dasharrow}
\newcommand{\divides}{\mid}
\newcommand{\iso}{\cong}
 \newcommand{\into}{\hookrightarrow}
 \newcommand{\onto}{\twoheadrightarrow}

\DeclareMathOperator{\divi}{div}
\DeclareMathOperator{\hcf}{hcf}
\DeclareMathOperator{\rank}{rank}
\DeclareMathOperator{\wt}{wt}
\DeclareMathOperator{\Cl}{Cl}
\DeclareMathOperator{\Grass}{Grass}
\DeclareMathOperator{\Pf}{Pf}
\DeclareMathOperator{\Pic}{Pic}
\DeclareMathOperator{\Proj}{Proj}
\DeclareMathOperator{\Spec}{Spec}
\DeclareMathOperator{\OGr}{OGr}
\DeclareMathOperator{\wGr}{wGr}

\newcommand{\FFbar}{\overline{\mathbb F}}

\newcommand{\Case}[1]{\paragraph{\sc Case #1}}
\newcommand{\rest}[1]{{}_{{\textstyle|}#1}}

 \theoremstyle{plain}
 \newtheorem{thm}[subsection]{Theorem}
 \newtheorem{thm-dfn}[subsection]{Theorem-Definition}

 \theoremstyle{definition}
 \newtheorem{dfn}[subsection]{Definition}
 \newtheorem{exa}[subsection]{Example}
 \newtheorem{exc}[subsection]{Exercise}

 \theoremstyle{remark}
 \newtheorem{rmk}[subsection]{Remark}
 \newtheorem{rmks}[subsection]{Remarks}

\numberwithin{equation}{subsection}
\numberwithin{figure}{subsection}

\title{Fano 3-folds, K3 surfaces and graded rings}
\author{Selma Alt{\i}nok}
\address{Selma Alt{\i}nok, Adnan Menderes University, Art and Science
Faculty, Department of Mathematics, Aydin 09010, Turkey}
\email{saltinok43@hotmail.com}
\author{Gavin Brown}
\address{Gavin Brown, Math Inst., Univ.\ of Warwick, Coventry CV4 7AL,
England}
\email{gavinb@maths.warwick.ac.uk}
\urladdr{www.maths.warwick.ac.uk/$\!\scriptstyle\sim$gavinb}
\author{Miles Reid}
\address{Miles Reid, Math Inst., Univ.\ of Warwick, Coventry CV4 7AL,
England}
\email{miles@maths.warwick.ac.uk}
\urladdr{www.maths.warwick.ac.uk/$\!\scriptstyle\sim$miles}
\date{Feb 2002}
 
\begin{document}

 \begin{abstract}
 Explicit birational geometry of 3-folds represents a second phase of
Mori theory, going beyond the foundational work of the 1980s. This
paper is a tutorial and colloquial introduction to the explicit
classification of Fano \hbox{3-folds} (also known by the older name
$\Q$-Fano 3-folds), a subject that we hope is nearing completion.
With the intention of remaining accessible to beginners in algebraic
geometry, we include examples of elementary calculations of graded
rings over curves and K3 surfaces. For us, K3 surfaces have at worst Du Val
singularities and are polarised by an ample Weil divisor (you might
prefer to call these $\Q$-K3 surfaces); they occur as the general
elephant of a Fano 3-fold, but are also interesting in their own right.
A second section of the paper runs briefly through the classical theory
of nonsingular Fano 3-folds and Mukai's extension to indecomposable
Gorenstein Fano \hbox{3-folds}. Ideas sketched out by Takagi at the
Singapore conference reduce the study of $\Q$-Fano 3-folds with $g\ge2$
(and a suitable assumption on the general elephant) to indecomposable
Gorenstein Fano 3-folds together with unprojection data.

Much of the information about the anticanonical ring of a Fano
\hbox{3-fold} or K3 surface is contained in its Hilbert series. The
Hilbert function is determined by orbifold Riemann--Roch (the
Lefschetz formula of Atiyah, Singer and Segal, see Reid \cite{YPG});
using this, we can treat the Hilbert series as a simple collation of
the genus and a basket of cyclic quotient singularities. Many hundreds
of families of K3s and Fano 3-folds are known, among them a large
number with $g\le0$, and Takagi's methods do not apply to these.
However, in many cases, the Hilbert series already gives firm
indications of how to construct the variety by biregular or birational
methods. A final section of the paper introduces the K3 database in
Magma, that manipulates these huge lists without effort.
 \end{abstract}

\maketitle
 \setcounter{tocdepth}{1}
 \tableofcontents

\section{Introduction}
\subsection{The graded ring $R(X,A)$ of a polarised variety}
Let $X$ be an irreducible projective variety over $\C$, and $A$ an
ample divisor on $X$ (see below for more explanation). For $n\ge0$, we
write
 \[
 H^0(X,nA) = H^0(X,\Oh_X(nA)) = \sL(nA) =
 \bigl\{ f \in \C(X) \bigm| \divi f + nA \ge 0 \bigr \}
 \]
 for the {\em Riemann--Roch space} (RR space) of $nA$. That is,
$H^0(X,nA)$ is the finite dimensional vector space of rational (or
meromorphic) functions $f\in\C(X)$ with divisor of poles $\le nA$. Our
basic construction is the graded ring
 \begin{equation}
 R(X,A) = \bigoplus_{n\ge0} H^0(X,nA),
 \label{eq!RX}
 \end{equation}
where the product is simply multiplication of rational functions
 \[
 H^0(X,nA) \times H^0(X,mA) \to H^0(X,(n+m)A) \quad\hbox{by}\quad
 (f,g) \mapsto fg.
 \]
This just says that if $f,g$ are rational functions with poles $\le
nA,mA$ then $fg$ has poles $\le (n+m)A$.

Our special interest is the case when the ring $R(X,A)$ can be
described by explicit generators and relations, for example, as a
polynomial ring or a polynomial ring divided by a principal ideal
(geometrically, a hyper\-surface). When this is possible, it
corresponds to embedding $X$ in projective space and determining the
defining equations of the image variety. It frequently happens in
higher dimensions that the generators $x_i$ of $R(X,A)$ have different
weights, so we have to work with weighted projective spaces ({\em
w.p.s.})\ and weighted homogeneous ideals. Dolgachev \cite{D} and
Fletcher \cite{Fl} are useful as general references on w.p.s.\ and their
complete intersections ({\em c.i.}).

 \begin{rmks}
 Although we intend to allow $X$ to have mild singularities, please think
of it in the first instance as nonsingular. The point of the RR spaces
$H^0(X,nA)=\sL(nA)$ is this: because $X$ is projective, the only globally
defined regular functions (holomorphic functions) on it are the constants. By
allowing poles along a divisor $nA$, we get a finite dimensional space of
rational (or meromorphic) functions $H^0(X,nA)$; as we discuss below, the RR
theorem predicts $\dim H^0(X,nA)$ in good cases. Choosing a basis
$x_0,x_1,\dots,x_k\in H^0(X,nA)$ defines a rational map
 \[
 \fie_{nA}\colon X \broken \PP^k \quad\hbox{by}\quad
 P \mapsto (x_0(P):\cdots:x_k(P));
 \]
there are straightforward criteria that determine when $\fie_{nA}$ is a
morphism $X\to\PP^k$ or an embedding $X\into\PP^k$. The divisor $nA$ is
{\em very ample} if it defines an embedding $\fie_{nA}\colon X\into \PP^k$,
and $A$ is {\em ample} if $nA$ is very ample for some $n$.
 \end{rmks}
 
 \subsection{The RR theorem for a curve}
 Please skip this stuff if you already know it. If not, we strongly advise
you to make an effort to {\em commit it to memory}, since it is one of the
central points of algebraic geometry. A nonsingular projective curve $C$ (or
compact Riemann surface) has a {\em genus} $g$ that can be defined in several
alternative ways:
 \begin{enumerate}
 \renewcommand{\labelenumi}{(\arabic{enumi})}
 \item $C$ is homeomorphic to the traditional picture of a surface with $g$
holes.
 \item $C$ has Euler number $e(C)=2-2g$.
 \item The tangent bundle to $C$ has degree $2-2g$, or equivalently, the
canonical divisor class $K_C$ has degree $2g-2$. Here $K_C$ is defined as the
class of the divisor of a rational (or meromorphic) differential.
 \item $\dim H^0(C,\Oh_C(K_C))=\dim \sL(K_C)=g$. In words, $g$ equals
 the dimension of the vector space of global regular (or holomorphic)
 differentials on $C$.
 \item etc.; for example, you can define $g$ as the integer appearing in the
RR theorem for curves (\ref{eq!RRC}).
 \end{enumerate}

 \begin{thm}[RR for curves]
 Let $D=\sum n_iP_i$ be a divisor on $C$. Then
 \begin{equation}
 \dim H^0(D)-\dim H^0(K_C-D)=1-g+\deg D,
 \label{eq!RRC}
 \end{equation}
where $\deg D=\sum n_i$. (Recall that $H^0(D)=\sL(D)$ and
$H^0(K_C-D)=\sL(K_C-D)$, and in this language the formula is
 \[
 \dim \sL(D)-\dim \sL(K_C-D)=1-g+\deg D.)
 \]

When $\deg D>2g-2$, we get $H^0(K_C-D)=0$, so (\ref{eq!RRC}) simplifies
to a formula for $\dim H^0(D)$ in terms of topological invariants of
$C$ and $D$.
 \end{thm}

 \begin{exa} \label{exa!g=3}
 As a baby case of our graded ring methods, we compute the canonical ring
$R(C,K_C)$ of a curve of genus~3, and conclude that $C$ is a plane quartic
curve under suitable extra assumptions.

Take $C$ to be a nonsingular projective curve of genus~3, with polarising
divisor $K_C$. The starting point is that RR gives $\dim H^0(nK_C)$:
 \begin{equation}
 \dim H^0(nK_C)\ =\ 
 \begin{cases}
 1 & \hbox{if } n=0, \\
 g & \hbox{if } n=1, \\
 (2n-1)(g-1) & \hbox{if } n\ge2.
 \end{cases}
 \label{eq!g=3}
 \end{equation}
The last number comes from (\ref{eq!RRC}), with the r.h.s.\ equal to
$1-g+n(2g-2)$.

 Now $\dim H^0(C,K_C)=3$. As we said above, choosing a basis $x_1,x_2,x_3$
of $H^0(C,K_C)$ defines a rational map $\fie_{K_C}\colon C\broken\PP^2$.
A traditional and easy argument based on RR applied to the divisors $K_C-P$
and $K_C-P-Q$ proves that $\fie_{K_C}$ is a morphism, either defining an
isomorphism $C\iso C_4\subset\PP^2$ of $C$ with a plane quartic curve $C_4$,
or a generically 2-to-1 cover $C\to Q_2\subset\PP^2$ over a conic.

We argue somewhat more algebraically on the ring
 \[
 R(C,K_C)=\bigoplus_{n\ge0}H^0(C,nK_C)
 \]
 using (\ref{eq!g=3}); in degree two, $\dim H^0(C,2K_C)=6$. But we already
know 6 elements of $H^0(C,2K_C)$, namely the 6 quadratic monomials
 \[
 S^2x_i=\bigl\{x_1^2,x_1x_2,x_2^2,x_1x_3,x_2x_3,x_3^2\bigr\}.
 \]
We assume that these monomials are linearly independent, so form a basis
of $H^0(C,2K_C)$. Likewise, in degree 3, $\dim H^0(C,3K_C)=10$, and there
are $10=\binom{5}{2}$ cubic monomials $S^3x_i$, so we assume that they
are linearly independent, and again form a basis. In degree 4, however,
we necessarily find a relation, since $\dim H^0(C,4K_C)=14$, but there are
$15=\binom{6}{2}$ quartic monomials $S^4x_i$.

The conclusion is the prediction that the graded ring $R(C,K_C)$ has the
simplest form $R(C,K_C)=\C[x_1,x_2,x_3]/(f_4)$; and the corresponding
map $\fie_{K_C}\colon C\to\PP^2$ is an embedding with image
$C_4$ given by $f_4(x_1,x_2,x_3)=0$.
 \end{exa}

 \begin{rmk}
 The assumption that the quadratic monomials $S^2x_i$ are linearly
independent is of course the hyperelliptic dichotomy. In the
non\-hyper\-elliptic case, the image $\fie_{K_C}(C)$ cannot be contained in a
curve of degree $\le3$ union a finite set, so that once the relation $f_4$
has been detected, one sees that the ring homomorphism
 \[
 \C[x_1,x_2,x_3]/(f_4) \to R(C,K_C)
 \]
 must be injective. Then it is surjective, because in degree $n$ the
quotient ring has dimension $\binom{n+4}{2}-\binom{n}{2}=\dim H^0(C,nK_C)$.

In the hyperelliptic case, the canonical ring is
 \[
 \C[x_1,x_2,x_3,y]/(Q_2,F_4),
 \]
where $Q_2(x)=0$ is the equation of the image conic, $y$ the new generator
in degree 2 required to compensate for this relation, and $F_4:y^2=f_4(x)$.
Then $C\to Q_2$ is the double cover ramified in the 8 points $Q_2=f_4=0$.

We can even consider a degenerating family with quadratic relation
$\la y=Q_2(x)$ that consists of a nonhyperelliptic curve if $\la\ne0$ and a
hyperelliptic curve if $\la=0$.
 \end{rmk}

 \begin{exa}
 Consider now a nonhyperelliptic curve $C$ of genus~6. Its canonical
embedding $\fie\colon C_{10}\into\PP^5$ has image of codimension~4. As
opposed to the hypersurface case, we do not have any surefire way of
predicting the equations of a variety of codimension~$\ge4$. This is a
major preoccupation of the rest of the paper. In this case, we happen to
be lucky, but it still involves a case division into Brill--Noether
special and Brill--Noether general:

\Case1 $C$ is trigonal or isomorphic to a plane quintic. Trigonal means
that $C$ has a linear system $g^1_3$; equivalently, $C$ can be
represented as a triple cover $C\to\PP^1$. It is well known that then
$\fie_{K_C}(C)$ is contained in a quartic surface scroll (typically,
$\PP^1\times\PP^1$ embedded by $\Oh(1,2)$); its equations are the 6
quadrics defining the scroll and 3 cubics defining $C$ inside the
scroll. The case of a plane quintic is similar, with $\fie_{K_C}(C)$
contained in the Veronese image of $\PP^2$, and defined by the 6
quadratic equations of the Veronese surface and 3 cubics defining $C$
in it.

\Case2 $C$ is not trigonal and not isomorphic to a plane quintic. Then
the canonical image $\fie_{K_C}(C)$ is a c.i.\ in $\Grass(2,5)$:
 \[
 C\ =\ \Grass(2,5) \cap Q_2 \cap H_1 \cap \cdots \cap H_4 \subset\PP^9.
 \]
 \end{exa}

 \begin{rmks} \label{rmk!Pf}
 Here $\Grass(2,5)$ is the Grassmann variety of 2-dimensional vector
subspaces of $\C^5$ in its Pl\"ucker embedding in
$\PP^9=\PP(\bigwedge^2\C^5)$. Equivalently, it is the determinantal
variety defined by the $4\times4$ diagonal Pfaffians $\Pf_{ij.kl}$ of
the $5\times5$ skew matrix
 \[
 \begin{pmatrix}
 x_{12} & x_{13} & x_{14} & x_{15} \\
 & x_{23} & x_{24} & x_{25} \\
 && x_{34} & x_{35} \\
 &&& x_{45} \\
 \end{pmatrix},
 \]
 where $\pm\Pf_{ij.kl}=x_{ij}x_{kl}-x_{ik}x_{jl}+x_{il}x_{jk}$ for
$\{i,j,k,l\} \subset\{1,2,3,4,5\}$. (We only write the upper triangular
elements. For the whole $5\times5$ matrix, just write zeros in the
5 diagonal entries and the skew elements $-x_{ji}$.)

Although Case~2 is a 19th century result, the direct statement and proof
is due to Mukai (\cite{Mu1}, Section~5), and is a first substantial case
of his general program. He proves that there exists a unique rank~2
stable vector bundle $E$ on $C$ with $\det E=K_C$ and $H^0(C,E)=\C^5$;
moreover, $E$ is generated by its $H^0$, and $\bigwedge^2 H^0(C,E)\onto
H^0(C,\bigwedge^2E)$. By the universal mapping property of the Grassmann
variety, $E$ and its sections define a morphism $\psi_E\colon
C\to\Grass(2,5)$, making $\fie_{K_C}$ a linear section of the Pl\"ucker
embedding:
 \[
 \renewcommand{\arraystretch}{1.6}
 \begin{array}{ccc}
 C & \into & \PP^5 \\
 \big\downarrow && \bigcap \\
 \Grass(2,5) &\into& \PP^9
 \end{array}
 \]
Since $g=6\ge2\cdot3$, standard Brill--Noether arguments give the
existence of a decomposition $K_C=\xi+\eta$ with divisors $\xi$ and
$\eta$ such that $H^0(C,\xi)=2$ and $H^0(C,\eta)=3$. RR (\ref{eq!RRC})
gives $\deg\xi=4$, and the case assumptions imply that $|\xi|$ is a free
$g^1_4$ and $|\eta|$ a free $g^2_6$. For every such decomposition
$K_C=\xi+\eta$, the bundle $E$ appears as the unique extension $\xi\into
E\onto\eta$ for which the boundary map $H^0(\eta)\to H^1(\xi)$ is zero.
 \end{rmks}

 \begin{exa} \label{exa!1/2}
As in Example~\ref{exa!g=3}, let $C$ be a nonsingular quartic curve;
however, instead of $K_C$, we choose a point $P\in C$ and consider the
{\em fractional divisor} $A=K_C+\half P$. The RR space of $nA$ is given
by the same formula:
 \[
 H^0(C,nA)=\sL(nA)=\bigl\{f \in\C(C) \bigm| \divi f+nA\ge0 \bigr\}.
 \]
The novel point is that a rational function cannot have a pole of
fractional order. Thus although the definition allows $f$ a pole of
order $n/2$ at $P$, in practice this restricts it to have pole of
order $\rd{\frac n2}$.

The first time this has any effect is when $n=2$, when
 \[
 H^0(C,2A)=H^0(C,2K_C+P)
 \]
 has dimension 1 bigger than $H^0(C,2K_C)$. Thus the ring $R(C,A)$
needs a new generator $y$ in degree~2. In
explicit terms, if $P=(1,0,0)\in C\subset\PP^2$ then the equation of
$C$ is $f_4=x_2A+x_3B=0$ where $A,B$ are cubics in $x_1,x_2,x_3$,
and $y=B/x_2=-A/x_3$ is a rational section of $\Oh_C(2)$ with pole
at $P$. It follows that
 \[
 R(C,A)=\C[x_1,x_2,x_3,y]/(yx_2-B,yx_3+A).
 \]
 \end{exa}

 \begin{rmks} \label{rmk!orb}
 The ring $\C[x_1,x_2,x_3,y]/(yx_2-B,yx_3+A)$ constructed in
Example~\ref{exa!1/2} corresponds to embedding $C\into\PP(1,1,1,2)$ as
a $(3,3)$ c.i. Here the w.p.s.\ $\PP(1,1,1,2)$ is the cone over the
Veronese surface, and $C$ passes through the cone point $(0,0,0,1)$ as a
nonsingular branch, with the equations $yx_2=B,yx_3=-A$ determining
$x_2,x_3$ as implicit functions.

 The fractional divisor $K_C+\frac{r-1}{r}P$ is the {\em orbifold
canonical class} of $C$, where $P$ is viewed as an orbifold point of
order~$r$. The ring $R(C,A)$ is again Gorenstein (see Watanabe
\cite{W}); we treat the case $K_C+\frac{2}{3}P$ in
Example~\ref{exa!2/3}. More generally, we could take any {\em basket}
of orbifold points $P_i$ of order $r_i$ on $C$ and consider the
orbifold canonical class $A=K_C+\sum\frac{r_i-1}{r_i}P_i$. The affine
cone $\Spec R(C,A)$ corresponds to a $\C^*$ fibration over $C$ with
isotropy $\bmu_{r_i}$ over $P_i$, the $\C^*$ analog of a Seifert fibre
space. The fact that $R(C,A)$ is a Gorenstein ring means that there is
a category of orbifold sheaves on $C$ (more precisely, modules over its
affine cone $\sC C=\Spec R(C,A)$) having a nice form of Serre
duality.

 A curve with an orbifold point is a substantial first case of the cyclic
quotient singularities needed in higher dimension. For example, our
construction of $C_{3,3}\subset\PP(1,1,1,2)$ extends in the obvious way
to a K3 surface $S_{3,3}\subset\PP(1^4,2)$, Fano 3-fold
$V_{3,3}\subset\PP(1^5,2)$, etc.
 \end{rmks}

 \section{Classic Fano 3-folds}
 The subject starts in the 1930s, when Fano studies projective 3-folds
having canonical curve sections. His definition is projective:
$V=V^3_{2g-2}\subset\PP^{g+1}$ should have canonical curves as its
codimension~2 linear sections:
 \[
 V\cap H_1\cap H_2 \ = \ C_{2g-2}\subset\PP^{g-1},
 \]
but this is more-or-less equivalent to assuming that $-K_V$ is very ample
(see below). Iskovskikh modernised Fano's treatment in the 1970s. His
starting point is a nonsingular 3-fold $V$ with $-K_V$ ample. He analyses
the rather few exceptions to $-K_V$ very ample, and corrects and reworks
Fano's classification. {From} the 1980s onwards Mukai discovers a new
interpretation of Fano and Iskovskikh's results in terms of linear
sections of projective homogeneous varieties, and generalises them to
indecomposable Gorenstein Fano 3-folds (see below for explanation).

 \begin{dfn}
 A {\em nonsingular Fano $3$-fold\/} is a nonsingular irreducible
projective variety of dimension $3$ with ample $-K_V$. Likewise, a {\em
Gorenstein Fano $3$-fold\/} is a normal irreducible projective $3$-fold
$V$ with at worst canonical Gorenstein singularities and $-K_V$ ample.
The {\em genus} of a Fano 3-fold is the integer $g$ defined by $\dim
H^0(V,-K_V)=g+2$. In the nonsingular case, we deduce that
$2g-2=({-}K_V)^3$.

 We say that $V$ is {\em prime} if $\Cl V=\Z\cdot({-}K_V)$ (here $\Cl
V$ is the Weil divisor class group); and $V$ is {\em indecomposable} if
there does not exist any decomposition $-K_V=A+B$ where $A,B$ are Weil
divisors such that $|A|,|B|$ are nontrivial linear systems.
 \end{dfn}

The point of these definitions is to exclude easy cases such as $\PP^3$
(for which $-K=\Oh(4)$) and quadric or cubic hypersurfaces $Q_2,
F_3\subset\PP^4$ (for which $-K=\Oh(3)$ or $\Oh(2)$ respectively), that
can be handled by simpler methods. The prime condition is Fano and
Iskovskikh's assumption that $V$ is factorial and has rank
$\rho(V)=\rank\Pic V=1$ and index~1; indecomposable is Mukai's
generalisation. Nonsingular Fano 3-folds with $\rho\ge2$ were treated
in detail by Mori and Mukai \cite{MM1}--\cite{MM2}.

 \begin{exa}\label{exa!V4}
 A quartic $3$-fold containing a plane $\Pi\subset V_4\subset\PP^4$ is
indecomposable. If $\Pi$ is the $x_3,x_4,x_5$ coordinate plane defined
by $x_1=x_2=0$ then $V$ is defined by an equation $x_1A+x_2B=0$, and
this is singular at the points $\Pi\cap(A=B=0)$ (in general 9 ordinary
double points).
 \end{exa}

 \subsection{General theory} \label{ssec!ele}
 An {\em elephant} of $V$ is a surface $S\in|{-}K_V|$ with at worst Du
Val singularities (rational double points). If $V$ is a Gorenstein Fano
3-fold, an elephant is known to exist by a theorem of Shokurov and Reid
\cite{R1}. It is a K3 surface: in fact $K_S=(K_V+S)\rest S=0$ and
$H^1(S,\Oh_S)=0$ by Kodaira vanishing. Morever, $|{-}K_V|\rest S$ is an
ample complete linear system of Cartier divisors on $S$.

 With a couple of exceptional cases that are easily classified and that
we pass over, it follows that $|{-}K_V|$ is very ample. Taking another
elephant $S'\in|{-}K_V|$ and setting $C=S\cap S'$ gives the diagram
 \[
 \renewcommand{\arraystretch}{1.3}
 \begin{array}{ccll}
 V&\into&\PP^{g+1} & \hbox{anticanonical embedding $\fie_{{-}K_V}$,}\\
 \bigcup&&\bigcup \\
 S&\into&\PP^g & \hbox{embedding by $|{-}K_V|\rest S$,} \\
 \bigcup&&\bigcup \\
 C&\into&\PP^{g-1} & \hbox{canonical embedding $\fie_{K_C}$.} 
 \end{array}
 \]
Special linear systems on $C$ interpreted in terms of geometric RR
give rise to geometric properties of $V$. The hyperelliptic case was
already passed over in what we said above; the case of $C$ trigonal
leads to $V$ a hypersurface in a scroll, which contradicts $V$
indecomposable for $g\ge5$. In the same way, other Brill--Noether
special linear systems on $C$ such as $g^2_5$ (when $C$ is a plane
quintic) usually contradict $V$ indecomposable.

 \subsection{Fano 3-folds, specific theory}\label{ssec:specific}
 None of the above is specific to 3-folds. Canonical curves and K3 surfaces
continue to exist for all $g$, so the main result for 3-folds is thus quite
remarkable.

 \begin{thm}[Fano, Iskovskikh, Mukai] \label{th!clF}
 Indecomposable Gorenstein Fano \hbox{$3$-folds} of genus $g$
exist if and only if $2\le g\le10$ or $g=12$. If $g\le5$, the
anticanonical ring $R(V,{-}K_V)$ is a hypersurface or complete
intersection in projective space or w.p.s. If $g=6,\dots,10$ then $V$ is
a complete intersection in a homogeneous projective variety. (There is
also an analogous structure result for
$g=12$.)
 \end{thm}

For example,
 \[
 \begin{array}{rcl}
 g=2 & \implies & V=Q_6\subset\PP(1,1,1,1,3), \\[6pt]
 g=4 & \implies & V=Q_2\cap F_3\subset\PP^5.
 \end{array}
 \]
As a typical case of the homogeneous projective varieties $G/P$,
 \[
 \begin{array}{rcl}
 g=7 & \implies & V=\Si\cap H_1\cap \cdots\cap H_7,
 \end{array}
 \]
where $\Si=\OGr^{10}(5,10)$ is the 10-dimensional spinor variety or
orthogonal Grassmann variety, that is, the subset of $\Grass(5,10)$
consisting of maximal isotropic vector subspaces of the standard inner
product $\left(\begin{smallmatrix} 0 & I \\ I & 0
\end{smallmatrix}\right)$ in its spinor embedding $\Si\subset\PP^{15}$
(see Mukai \cite{Mu2}).

As an illustration, we sketch the proof in the style of Mukai that
there does not exist any Fano 3-fold $V=V_{20}\subset\PP^{12}$ of genus
11. Write $\si\colon V_1\to V$ for the blowup of a point $P\in V$, and
$E\subset V_1$ for the exceptional divisor, with $E\iso\PP^2$ and
$\Oh_E(-E)\iso\Oh_{\PP^2}(1)$. Then $K_{V_1}=\si^*K_V+2E$. It follows
that the anticanonical linear system $|{-}K_{V_1}|=|\si^*({-}K_V)-2E|$
is the birational transform of the linear system
$|m_P^2\cdot\Oh_V({-}K_V)|$ of hyperplane sections of $V$ containing
the tangent plane $T_PV$. We assume that $P$ does not lie on any line
of $V$. Then since $V$ is an intersection of quadrics, we deduce that
$\fie_{{-}K_{V_1}}\colon V_1\to\PP^8$ is a morphism, birational because
$V$ is indecomposable, with image $V'$ a Gorenstein Fano 3-fold of
genus~7. Now on the one hand, $V'$ is indecomposable, so is a linear
section of $\OGr^{10}(5,10)$ by the result for $g=7$. On the other
hand, it contains the image of $E$, a Veronese surface.

Finally, arguing on the geometry of $\OGr^{10}(5,10)\subset\PP^{15}$ and
on its universal bundle, we can deduce that the only Veronese surfaces $E$
it contains span a \hbox{4-plane} $\PP^4\subset\OGr$; this contradicts the
construction of $V'$ as a linear section of $\OGr$ containing $E$.

Iskovskikh's proof (deriving from Fano) proceeds by the projection from
a general line $L\subset V$; it involves proving that lines exist, and
various generality statements about its projection and double
projection. Mukai's proof is thus a considerable simplification even in
the nonsingular case.

 \section{$\Q$-Fano 3-folds}
 \subsection{The Mori category}
 The minimal model program for 3-folds (usually called MMP or Mori theory)
was developed by Mori and others from the late 1970s. As in the theory
of surfaces, to get one step closer to a minimal model, one makes a
birational contraction, for example, contracting a copy of $\PP^2$ on
which the canonical class is negative. However, these contractions lead to
singularities, so that Mori theory only works in a suitable category of
singular varieties. This leads to the following definitions:
 \begin{dfn}
 The {\em Mori category} consists of projective varieties with at worst
\hbox{$\Q$-factorial} terminal singularities, where
 \begin{enumerate}
 \renewcommand{\labelenumi}{(\alph{enumi})}
 \item A variety $V$ has {\em terminal singularities} if it is normal,
$rK_V$ is a Cartier divisor for some $r>0$, and any resolution of
singularities $f\colon Y\to V$ with divisorial exceptional locus
$\bigcup E_i$ satisfies
 \[
 K_Y = f^*K_V + \sum a_iE_i \quad\hbox{with all $a_i>0$.}
 \]
For example, the cone over the Veronese surface (that is, the quotient
singularity $\half(1,1,1)$) is terminal: it is resolved by the blowup
$Y\to V$ of the origin, which introduces the exceptional divisor
$E\iso\PP^2$ with $\Oh_E(-E)\iso\Oh_{\PP^2}(2)$ and $K_Y{}\rest E\iso
\Oh_{\PP^2}(-1)$, so that $K_Y=f^*K_X+\half E$.

 \item $V$ is {\em $\Q$-factorial\/} if every Weil divisor $D$ on $V$
has a multiple $rD$ that is Cartier. You should think of this as an
analog in algebraic geometry of the condition that $V$ is a
$\Q$-homology manifold: a codimension~1 sub\-variety has a dual
cohomology class in $H^2(V,\Q)$. For example, the quartic hypersurface
of Example~\ref{exa!V4} is not $\Q$-factorial, since it has double
points on $\Pi$ at which no multiple of $\Pi$ is Cartier.
 \end{enumerate}
 \end{dfn}

 \subsection{Terminal singularities}
 Terminal singularities are necessary for Mori theory: the Mori category
is closed under contractions and flips, and is the smallest such category.

 Terminal singularities were classified by Mori and Reid. They are made
up of the following ingredients (for more details, see Reid \cite{YPG}):
 \begin{enumerate}
 \renewcommand{\labelenumi}{(\arabic{enumi})}
 \item Mild isolated hypersurface singularities, for example, double
points such as $(xy=f(z,t))\subset\C^4_{x,y,z,t}$. While this is a
reasonably concrete class of singularities, already the special case
$xy=f(z,t)$ is fairly infinite, containing all isolated plane curve
singularities.
 \item The cyclic quotient singularities $\recip{r}(1,a,r-a)$. The
notation means the quotient $V=\C^3/(\Z/r)$, where the cyclic group
$\Z/r$ acts by
 \[
 (x,y,z)\mapsto (\ep x,\ep^ay,\ep^{r-a}z).
 \]
 Setting $x=0$ leads to the cyclic quotient singularity
$\recip{r}(a,r-a)$, which is the Du Val singularity A$_{r-1}$
 \[
 S:(uv=w^r)\subset \C^3_{u,v,w}.
 \]
Thus $(x=0)$ defines a local elephant $S\in|{-}K_V|$.
 \item Some combinations of the above two types, typically the
main Type~A family $\bigl(xy=f(z^r,t)\bigr)/\recip{r}(a,r-a,1,0)$.
 \end{enumerate}

 \begin{rmk} \label{rmk!a?b}
 Whenever we say that a point $P\in V$ is a terminal quotient
singularity of type $\recip{r}(1,a,r-a)$, we always assume that $0<a<r$
and $a$ is coprime to $r$. The local class group of $P\in V$ is $\Z/r$,
with chosen generator $\Oh_V(-K_V)=\Oh(1)$, a {\em polarisation} of
the singularity. The generators of the graded ring $R(V,A)$ that
serve as orbifold local coordinates then have weights $1,a,r-a$ mod
$r$; compare Corti, Pukhlikov and Reid \cite{CPR}, 3.4.6.

The quotient singularity $\recip{r}(a,r-a)$ is the Du Val singularity
A$_{r-1}$, so of course does not depend on $a$ up to analytic
isomorphism; but we work here with K3 surfaces with a chosen
polarisation $\Oh(1)=\Oh(D)$, and the weights mod $r$ of the orbifold
local coordinates are uniquely determined.

Write $b$ for the inverse of $a$ mod $r$ (that is $ab\equiv1$ mod $r$).
Choosing $\ep'=\ep^a$ as a new basis for the group of $r$th roots of
unity would put the Du Val singularity in the standard A$_{r-1}$ form
$\recip{r}(1,r-1)$ and the terminal 3-fold point in the form
$\recip{r}(b,1,r-1)$; and $\Oh_S(D)$ is locally isomorphic to the
$\ep'{}^b$ eigensheaf. The birational transform of the general divisor
$D$ meets the $b$th curve in the resolution (see
Figure~\ref{fig!star}). The quantity $b$ also appears in the orbifold
RR formulas of \cite{YPG}, Theorem~10.2 (see Theorem~\ref{th!Hi}
below), for essentially the same reason. Muddling up $a$ and $b$ in
calculations is a common error, but you find out soon enough when your
plurigenera turn out to be fractional; see Exercise~\ref{exc!ele}, (3)
for a practical instance.
 \end{rmk}

 \begin{dfn}
 A {\em Fano $3$-fold} is a variety $V$ for which $-K_V$ is ample. We
usually add conditions to this. A {\em Mori Fano $3$-fold} $V$ is a Fano
3-fold in the Mori category and with $\rank\Pic V=1$. As one of the
possible end products of a MMP, these 3-folds are among the basic
building blocks of Mori theory. $V$ is {\em prime} if it is in the Mori
category and $\Cl V=\Z\cdot(-K_V)$.

 As in the Gorenstein case, the {\em anticanonical ring} of a Fano
3-fold $V$ is the graded ring
 \begin{equation}
 R(V,-K_V)=\bigoplus_{n\ge0} H^0(V,-nK_V)).
 \label{eq!RV}
 \end{equation}
 An {\em elephant} of $V$ is a general divisor $S\in|{-}K_V|$, just as
in \ref{ssec!ele}. If $S$ exists, it is no longer a Cartier divisor at
singularities of $V$ of index $r>1$. The graded rings of $V$ and $S$
are related by
 \begin{equation}
 R(S,A) = R(V,-K_V)/(x_0),
 \label{eq!RS}
 \end{equation}
where $x_0\in H^0(V,{-}K_V)$ is the equation of $S\in|{-}K_V|$. That is,
$R(S,A)$ is a quotient of $R(V,-K_V)$ by a principal ideal generated by
an element of degree~1. When two rings are related in this way, many
basic algebraic properties of one can be inferred from the other;
compare Exercise~\ref{exc!ele} and \ref{sssec!F3}. This is called the
{\em hyperplane section principle.} The most useful case is when there
is an elephant $S$ with at worst Du Val singularities; then $S$ is a K3
surface (this is proved as in \ref{ssec!ele}) polarised by the Weil
divisor $A=-K_V{}\rest S$. Then via the graded rings, many properties
of the Fano 3-fold $V$ can be read from those of $S$.
 \end{dfn}

 \begin{rmk}
 As we see in \ref{sec!g<0} below, there are cases when
$|{-}K_V|=\emptyset$, that is, no elephant exists; or when
$h^0(-K_V)=1$, and the single surface $S=|{-}K_V|$ happens to have
slightly worse than Du Val singularities, so there is no K3 elephant.
Nevertheless, the motivation arising from K3 surfaces is our main
guiding principle for the study of Fano 3-folds, even in cases such as
these when it is logically inapplicable. Compare Exercise~\ref{exc!ele}.
 \end{rmk}

 \begin{exa} \label{exa!2/3}
 Let $x_1,\dots,x_5,y$ be coordinates on $\PP(1^5,2)$. We start from a
$(3,3)$ c.i.\ $V_{3,3}\subset\PP^5(1^5,2)$ as described at the end of
Remark~\ref{rmk!orb}, and assume that $V$ contains the weighted
projective plane
 \[
 \Pi=\PP(1,1,2) \quad\hbox{given by}\quad (x_1=x_2=x_3=0),
 \]
and is otherwise general. Since $\Pi$ is a c.i., the two cubic equations
of $V$ are of the form
 \begin{equation}
 V: \begin{pmatrix}
 a_1&a_2&a_3 \\
 b_1&b_2&b_3 
 \end{pmatrix}
 \begin{pmatrix} x_1 \\ x_2 \\ x_3
 \end{pmatrix} =0,
 \label{eq!3}
 \end{equation}
where $a_i,b_i$ are forms of degree 2 in $x_1,\dots,x_5,y$. We show
that
$\Pi$ can be contracted to a quotient singularity of type
$\recip{3}(1,1,2)$ by a rational map $V\broken W$, where
$W\subset\PP(1^5,2,3)$ is a Mori Fano 3-fold. The new coordinate $z$
will be a rational section of $\Oh_V(-3K_V)$ with pole along $\Pi$.
Thinking of (\ref{eq!3}) as 2 simultaneous equations for $x_1,x_2,x_3$
with coefficients $a_i,b_i$ and solving them by Cramer's rule gives:
 \begin{equation}
 \begin{matrix}
 zx_1&=&a_2b_3-a_3b_2, \\
 zx_2&=&a_3b_1-a_1b_3, \\
 zx_3&=&a_1b_2-a_2b_1. \\
 \end{matrix}
 \label{eq!4}
 \end{equation}
The ``constant of proportionality''
 \[
 z=\frac{a_2b_3-a_2b_3}{x_1}=
 \frac{a_3b_1-a_3b_1}{x_2}=\frac{a_1b_2-a_1b_2}{x_3}
 \]
is a well defined rational form of degree~3 on $V$ with pole along $\Pi$.

Now consider the rational map
 \[
 \fie\colon V\broken W\subset\PP(1^5,2,3) \quad
 \hbox{given by} \quad (x_1,\dots,x_5,y,z);
 \]
 the image $W$ is given by the 5 equations (\ref{eq!3}--\ref{eq!4}).
$\fie$ is an
morphism wherever $\Pi$ is a Cartier divisor on $V$, contracts $\Pi$
to the point $(0,\dots,0,1)$, and is an iso\-morphism outside $\Pi$. The
point $(0,\dots,0,1)$ is the $\recip{3}(1,1,2)$ singularity
$\C^3_{x_4,x_5,y}/(\Z/3)$, because when $z=1$, (\ref{eq!4}) gives
$x_1,x_2,x_3$ as implicit functions. $V$ is quasi-smooth at points of
$\Pi$ where the matrix in (\ref{eq!3}) has rank~2, and $\Pi$ is a
Cartier divisor there. We can assume that this matrix has rank $\ge1$
everywhere on $\Pi$, hence the blowup of $\Pi$ is a small morphism. It
makes $\Pi$ a Cartier divisor so that $\fie$ is a morphism on the
blowup.

Notice that the graded ring $R(V,-K_V+\recip{3}\Pi)=R(W,-K_W)$ is generated
by $(x_1,\dots,x_5,y,z)$ with relations (\ref{eq!3}--\ref{eq!4}). These can
be put together as the five $4\times4$ Pfaffians of the skew matrix
 \[
 \begin{pmatrix}
 z&a_1&a_2&a_3 \\
 &b_1&b_2&b_3 \\
 &&x_3&-x_2 \\
 &&&x_1 \\
 \end{pmatrix} \quad\hbox{of degrees}\quad
 \begin{pmatrix}
 3&2&2&2 \\
 &2&2&2 \\
 &&1&1 \\
 &&&1 \\
 \end{pmatrix}.
 \]
See Remark~\ref{rmk!Pf} for our conventions on Pfaffians.

At the end of Remark~\ref{rmk!orb}, we noted that the Fano 3-fold
$V_{3,3}\subset\PP(1^5,2)$ has as its linear section the graded ring
$R(C,K_C+\half P)$ corresponding to the orbifold canonical class of a
curve of genus~3 together with an orbifold point $P$ of order $r=2$. The
construction of this example restricted to $C$ is the graded ring
$R(C,K_C+\frac23P)$ corresponding to the orbifold canonical class
$K_C+\frac23P$ (compare Example~\ref{exa!1/2}).
 \end{exa}

 \section{Numerical data and Hilbert series} \label{sec!Hi}
 \subsection{The aim}
 Several hundred families of K3 surfaces $S$ and Mori Fano \hbox{3-folds}
$V$ are known. We now develop methods to guarantee gratification from
this cornu\-copia. The cases treated in Examples~\ref{exa!1/2}
and~\ref{exa!2/3} are pretty tame, and more complicated things like
Example~\ref{exa!44} below are more typical. Since we study a polarised
variety $X,A$ in terms of the graded ring $R(X,A)$ of (\ref{eq!RX}),
the natural numerical invariants of $X$ to take is the {\em list of
ingredients} that go into its Hilbert series $P_X(t)$. We explain this
in Theorem-Definition~\ref{th!Hi}, but first we work through an example.

 \begin{exa} \label{exa!44}
 Consider the weighted hypersurfaces:
 \begin{equation}
 S_{44} \subset \PP(4,5,13,22) \quad\hbox{and}\quad
 V_{44} \subset \PP(1,4,5,13,22);
 \end{equation}
 this is Fletcher \cite{Fl}, List~13.3, No.~91, Reid1(91) in the Magma
database. We sketch how to calculate with these hypersurfaces. For
more details, see \cite{Fl}, Section~13 (compare Dolgachev \cite{D}).
Write $x_0,\dots,x_4$ or $x,y,z,t,u$ for variables of weights
$a_0,\dots,a_4=1,4,5,13,22$, and let $f_{44}(x,y,z,t,u)$ be a general
weighted polynomial of degree~44 in them. Set
 \begin{itemize}
 \renewcommand{\labelitemi}{}
 \item $R[V]=\C[x,y,z,t,u]/\bigl(f_{44}(x,y,z,t,u)\bigr)$,
 \item $\sC V=\Spec R[V]:\bigl(f_{44}=0\bigr) \subset \C^5$, and
 \item $V=\Proj R[V]:\bigl(f_{44}=0\bigr) \subset \PP(1,4,5,13,22)$.
 \end{itemize}
Then $\sC V$ is the affine cone over $V$, and $V=\Proj R[V]$ the quotient
of $\sC V$ by the $\C^*$ action $x_i\mapsto \la^{a_i}x_i$ for $\la\in\C^*$
with the given weights $a_i$. By a combination of Bertini's theorem at
general points and explicit calculations at the coordinate strata, one
checks that $\sC V$ is nonsingular outside the origin for general
$f_{44}$. We say that $V$ is {\em quasi-smooth}. Then the only
singularities of $V\subset\PP(1,4,5,13,22)$ arise from the fixed
points of the $\C^*$ action; these are the points of the $x_i$-axis,
fixed by the cyclic group $\bmu_{a_i}\iso\Z/(a_i)$ of $a_i$th roots of
1, and the points of the $y,u$-plane, which are fixed by
$\bmu_2=\{\pm1\}$, because $2=\hcf(4,22)$. Thus in general $V$ has the
following cyclic quotient singularities:
 \begin{equation}
 \begin{cases}
 \half(1,1,1) & \hbox{at $(f=0)$ on the $yu$ line (one point);} \\
 \recip{5}(1,3,2) & \hbox{at the $z$ point $(0,0,1,0,0)$;} \\
 \recip{13}(1,4,9) & \hbox{at the $t$ point $(0,0,0,1,0)$.} \\
 \end{cases}
 \label{eq!44}
 \end{equation}
For example, along the $t$ axis of $\sC V$, we can assume that $f_{44}$
has the monomial $zt^3$, which means that $\frac{\partial f}{\partial
z}\ne0$ when $t\ne0$, so that $z$ is an implicit function of the other
variables. It follows that $\sC V$ is nonsingular there and that $V$ has
a quotient singularity of type $\recip{13}(1,4,22)=\recip{13}(1,4,9)$
at $(0,0,0,1,0)$. At each singular point, $\Oh_V(1)=\Oh_V(-K_V)$ is a
generator of the local class group.
 \end{exa}

 \begin{dfn}[Hilbert series]
Let $R=\bigoplus_{n\ge0} R_n$ be a graded ring. We assume that $R$ is
generated by finitely many elements $x_i$ of positive degree over $R_0=\C$.
Its Hilbert series $P(t)$ is defined by setting
 \[
 P_n=\dim_\C R_n \quad\hbox{and} \quad P(t)=\sum_{n\ge0} P_nt^n.
 \]
In cases of interest, $P_n$ is given by a formula of orbifold RR type
(see \cite{YPG}, Chapter~3), possibly with some corrections for low
values of $n$ when cohomology is still present.
 \end{dfn}

 \begin{exc}
 Use (\ref{eq!g=3}) to show that the canonical ring $R(C,K_C)$ of a curve
$C$ of genus $g$ has Hilbert series
 \[
 P(t) \ =\ \frac{1+(g-2)t+(g-2)t^2+t^3}{(1-t)^2}\,.
 \]

 Now write $A=K_C+\sum\frac{r-1}{r}P$ for the orbifold canonical class
of Remark~\ref{rmk!orb}; the sum takes place over a basket
$\sB=\{P,\frac{r-1}{r}\}$ of orbifold points $P$ of order $r$. In
degree $n$, only the rounded down integral divisor
$\rd{nA}=nK_C+\sum\rd{\frac{n(r-1)}r}P$ moves. Deduce that the Hilbert
series of the orbifold canonical ring $R(C,A)$ is given by
 \begin{align}
 P(t) & = \sum_\sB h^0\Bigl(nK_C+\sum\rd{\frac{n(r-1)}r}P\Bigr)t^n
\notag \\
 & = \frac{1+(g-2)t+(g-2)t^2+t^3}{(1-t)^2}+
 \sum_\sB\frac{t(t+\cdots+t^{r-1})}{(1-t)(1-t^{r})} 
 \label{eq!orbC1} \\
 & = \frac{1-(g-1)t-t^2}{1-t}+\frac{t}{(1-t)^2}\deg A
 -\sum_\sB\frac{1}{(1-t^r)}\sum_{i=1}^{r-1}
\textstyle{\frac{r-i}{r}}t^i. \notag
 % \label{eq!orbC2}
 \end{align}
The first expression (\ref{eq!orbC1}) shows how much the integral
divisor $\sum\rd{\frac{n(r-1)P}r}$ contributes to each $H^0(C,nA)$. The
second expression is a transparent case of an orbifold RR formula: we
write $\deg A=\frac1r\deg(rA)$ where $rA$ is a Cartier divisor ($A$
itself does not make sense as a sheaf); in each term, the effect of the
denominator $1-t^r$ is to multiply by $1+t^r+t^{2r}+\cdots$, so that
the fractional part $\{\frac{i(r-1)}{r}\}=\frac{r-i}{r}$ discarded in
the rounddown is repeated periodically with period $r$.

This formula for orbifold curves can serve as a model for the formulas
of \cite{YPG}, Chapter~3. In fact, a result of Becky Leng's thesis
\cite{Leng} derives the orbifold RR formula for 3-folds (\cite{CPR},
Theorem~10.2) via a reduction to the curve case.
 \end{exc}

 \subsection{Hilbert series for K3 surfaces and Fano 3-folds}
 The analogous formulas for the Hilbert series of the graded ring over a
K3 surface $S$ or a Fano 3-fold $V$ are contained in Alt{\i}nok \cite{A}
and \cite{A1}. The proofs are based on the results of \cite{YPG},
Chapter~3.

 \begin{thm-dfn}\label{th!Hi}
 Let $S$ be a K3 surface with Du Val singularities, and $D$ a Weil
divisor. The Hilbert series $P_S(t)=\sum_{n\ge0}h^0(S,nD)t^n$ is given
by
 \begin{equation}
 P_S(t)\ =\ \frac{1+t}{1-t}+\frac{t(1+t)}{(1-t)^3}\,\frac{D^2}2
 -\sum_\sB \recip{(1-t^r)}\sum\limits_{i=1}^{r-1}
\frac{\overline{bi}(r-\overline{bi})t^i}{2r}\,,
 \label{eq!K3-s}
 \end{equation}
where the sum takes place over a {\em basket} $\sB=\{\recip{r}(a,-a)\}$
of cyclic quotient singularities. In each term, $b$ is the inverse
of $a$ mod $r$, as in Remark~\ref{rmk!a?b}, and $\overline{bi}$ denotes
the minimal nonnegative residue mod $r$.

 We define the {\em genus} $g=g(S,D)$ by the formula:
 \[
 P_1(S,D)=h^0(\Oh_S(D))=g+1.
 \]
 The {\em numerical data} of $S,D$ are the genus $g$ and the basket
$\sB=\{\recip{r}(a,-a)\}$. Calculating $P_1(S,D)$ as the coefficient of
$t$ in $(\ref{eq!K3-s})$ gives a formula for $D^2$ in terms of $g$ and
$\sB$:
 \[
 D^2 = 2g-2 + \sum_\sB \frac{b(r-b)}{r}\,.
 \]

 Let $V$ be a Fano $3$-fold and write $A=-K_V$. Then its anticanonical
ring has Hilbert series
 \begin{align}
 P_V(t)&\ =\ \sum_{n\ge0}h^0(V,nA)t^n \notag \\
 &\ =\ \frac{1+t}{(1-t)^2}+\frac{t(1+t)}{(1-t)^4}\,\frac{A^3}2
 -\sum_\sB \recip{(1-t)(1-t^r)}
 \sum\limits_{i=1}^{r-1}
 \frac{\overline{bi}(r-\overline{bi})t^i}{2r}\,.
 \label{eq!Fano_3-f}
 \end{align}
where the sum takes place over a basket $\sB=\{\recip{r}(1,a,r-a)\}$ of
terminal quotient singularities. We define the {\em genus} $g=g(V)$ by
 \[
 P_1(V,A)=h^0(V,A)=h^0(\Oh_V(-K_V))=g+2.
 \]
 Then, as in the K3 case, the {\em numerical data} of $V$ are the genus
$g$ and the basket $\sB$; calculating $P_1$ from $(\ref{eq!Fano_3-f})$
gives essentially the same formula as above:
 \[
 A^3 = 2g-2 + \sum_\sB \frac{b(r-b)}{r}\,.
 \]
 \end{thm-dfn}

 \begin{exc} \label{exc!ele}
 \begin{enumerate}
 \renewcommand{\labelenumi}{(\arabic{enumi})}
 \item If a Fano 3-fold $V$ has a K3 elephant $S\in|{-}K_V|$, the
numerical data of $V$ and $S$ are related in the obvious way, and
(\ref{eq!K3-s}) and (\ref{eq!Fano_3-f}) only differ by the extra factor
$(1-t)$ in the denominator. Use the hyperplane section principle
(\ref{eq!RS}) to derive (\ref{eq!Fano_3-f}) from (\ref{eq!K3-s}).

 \item In the same way, if $S,C$ is a polarised K3 surface having only
singularities of type $\recip{r}(1,r-1)$ (with $C$ passing through each
singular point as a nonsingular curve in the $1$ eigenspace), derive
(\ref{eq!K3-s}) from the orbifold canonical curve formula
(\ref{eq!orbC1}). Singularities of type $\recip{r}(1,r-1)$ account for
the large majority of all cyclic quotient singularities in K3 baskets.

 \item The varieties in Example~\ref{exa!44} have $g=0$ and the basket
(\ref{eq!44}). Check that
 \[
 D^2=-4+\half+\frac{2\cdot3}6+\frac{3\cdot10}{13}=
 \frac{44}{4\cdot5\cdot13\cdot22}=\recip{130}
 \]
(note the shift from $a=4$ in $\recip{13}(4,9)$ to its inverse $b=10$
in the fractional contribution $\frac{3\cdot10}{13}$). Check that
$P_S(t)$ given by (\ref{eq!K3-s}) satisfies
 \[
 (1-t^4)(1-t^5)(1-t^{13})(1-t^{22})P_S(t)=1-t^{44}.
 \]
 \end{enumerate}
 \end{exc}

 \subsection{Numerical data of K3s and Fano 3-folds}\label{sec!g<0}
 \subsubsection{Invariants and inequalities}\label{sssec!ineq}
 As we saw in Theorem-Definition~\ref{th!Hi}, the numerical data of a
K3 surface or a Fano 3-fold consist of an integer genus $g$ and a
basket $\sB$ of fractional expressions. For a K3 surface in
characteristic~0 we have
 \begin{equation}
 g\ge-1, \quad 2g-2+\sum\frac{b(r-b)}r>0\quad\hbox{and}\quad
 \sum_\sB (r-1)\le 19.
 \label{ineq!K3}
 \end{equation}
 The last inequality comes from the fact that the singularities of the
basket contribute $\sum(r-1)$ exceptional $-2$-curves to $\Pic S$, that
form a negative definite set; in characteristic $p$, the result would
be $\sum(r-1)\le21$.

There are approximately 6,640 possible baskets satisfying
(\ref{ineq!K3}), and countably many possible values of $g$ for each
basket. A small number of extreme cases are excluded because the
polarisation and the basket cannot fit inside the K3 lattice $\LK$ (see
\ref{sssec!K3L} below); for example, a K3 surface cannot contain
$\ge17$ disjoint \hbox{$-2$-curves}. However, apart from these, all other $g$
and $\sB$ really occur. Most give rise to graded rings of high
codimension, and our methods based on graded ring are no longer
appropriate for studying them; the same applies to canonical curves or
K3 surfaces of genus $g\gg0$, for which explicit equations and graded
ring methods give little information, and other ideas are needed, such
as the period space and Torelli for K3 surfaces.

For a Fano 3-fold $V$ we have
 \begin{equation}
 g\ge-2, \quad 2g-2+\sum\frac{b(r-b)}r>0\quad\hbox{and}\quad
 \sum_\sB \left(r-\recip{r}\right)<24.
 \label{ineq!Fa}
 \end{equation}
The last inequality comes from \cite{YPG}, Corollary~10.3, the orbifold
analog of the RR formula $\recip{24}c_1c_2=\chi(\Oh_V)=1$; compare also
Kawamata \cite{Ka1}, 2.2 and \cite{Ka2}. In fact the argument of
\cite{Ka2} using Bogomolov stability gives
 \[
 \sum \sB\left(1-\recip{r}\right)<24-8(-K_V)^3,
 \]
a slightly stronger inequality.

Note that (\ref{ineq!Fa}) is stronger than (\ref{ineq!K3}) if the number
of singularities is large; for example, (\ref{ineq!Fa}) allows only $15$
singularities. On the other hand, (\ref{ineq!Fa}) is weaker if there are
few singularities -- for example, it would allow a cyclic singularity of
index~24, or a singularity of index 22 plus one of index 2. Since we
expect there to be rather few families of Fano 3-folds with $h^0(-K_V)=0$,
and those with $h^0(-K_V)>0$ to have a K3 elephant after deformation, it
seems unlikely that any of these extra Fano 3-fold cases really occur.

 \subsubsection{Negative genus}\label{ssec!uh}
 On a K3 surface, $g=-1$ means that the polarising divisor $D$ is
ineffective; or we can say that a curve of genus $-1$ does not exist
(pretty reasonable, uh?). A Fano 3-fold with $g=-2$ has
$|{-}K_V|=\emptyset$; we can think of $g=-2$ as saying that two things
do not exist: $V$ does not have a canonical curve section $S_1\cap S_2$,
and moreover, it does not even have an elephant $S\in|{-}K_V|$. Only
three or four families of Fano 3-folds with $h^0(-K_V)=0$ are known;
the simplest of these, and the first to be discovered, is the
codimension~2 weighted c.i.\ $V_{12,14}\subset \PP(2,3,4,5,6,7)$ due to
Fletcher (see \cite{Fl}, List~16.7, No.~60 and compare Reid \cite{Ki},
Example~9.14).

 \subsubsection{$\Q$-smoothing and the general elephant}
\label{sssec!GE}
 A Mori Fano 3-fold $V$ is allowed to have general terminal
singularities. Or it might happen that $|{-}K_V|\ne\emptyset$, but every
$S\in|{-}K_V|$ has an essential singularity (worse than Du Val); this
is not very surprising if $P_1(V)=1$. For example, one can construct a
surface $S_{14}\subset\PP(2,2,3,7)$ having an elliptic singularity,
but contained in a quasi-smooth $V_{14}\subset\PP(1,2,2,3,7)$.
(Explicit construction: start from the scroll $\FF_3$, with fibre $A$
and negative section $B$; then $S_{14}$ is the double cover ramified in
$B+C_1+C_2$ where $C_1\in|3A+B|$ is the general hyperplane section of
the cone $\FFbar_3$, and $C_2\in|7A+2B|$ has a tacnode on $C_1$, giving
$C_1+C_2$ an infinitely near triple point.)

It is known that a terminal 3-fold singularity has a {\em
$\Q$-smoothing}, that is, a small deformation with only cyclic
quotient terminal singularities (see \cite{YPG}, 6.4 for explicit
equations). This is the idea behind the basket appearing in
(\ref{eq!Fano_3-f}), and one step in its proof (\cite{YPG}, proof of
Theorem~10.2). It seems reasonable to conjecture that every prime Fano
3-fold $V$ also has a $\Q$-smoothing; compare Namikawa \cite{N} and
Minagawa \cite{Mi1}, \cite{Mi2} for partial results. One might hope to
prove this in terms of deformation theory and Hodge theoretic
consequences of $-K_V$ ample, by analogy with work on smoothing
Calabi--Yau 3-folds (see Gross \cite{G} and Namikawa and Steenbrink
\cite{NS}). Moreover, if $g\ge-1$, we also conjecture that $V$
together with its elephant $S\in|{-}K_V|$ has a {\em simultaneous
$\Q$-smoothing}, that is, a small deformation such that the pair
$S\subset V$ has only standard cyclic singularities
$\recip{r}(a,r-a)\subset\recip{r}(1,a,r-a)$. If in addition $g\ge0$ and
the basket of $V$ consists only of $\recip{r}(1,1,r-1)$, then after a
small deformation, $V$ would have a curve section $C=S_1\cap S_2$ that
is an orbifold canonical curve as in Remark~\ref{rmk!orb}.

 \subsubsection{Numerical data and the K3 lattice}\label{sssec!K3L}
 The numerical data of a polarised K3 surface $S$ corresponds to a
sublattice
 \[
 L\Bigl(g,\sB=\Bigl\{\recip{r}(a,r-a)\Bigr\}\Bigr)\subset\LK
 \]
 of the standard K3 lattice. We assume that, after a
\hbox{$\Q$-smoothing}, the singularities of $S$ are exactly the cyclic
quotient singularities of the basket. Let $f\colon T\to S$ be the
minimal resolution of singularities, and consider the sublattice
$L\subset\Pic T$ generated by the exceptional $-2$-curves of $f$
together with $f^*(ND)$, where $N$ is the global index of $D$. The
lattice $L=L(g,\sB)$ is based by the {\em quasistellar graph}
$\Ga(g,\sB)$ of Figure~\ref{fig!star}.
 \begin{figure}[ht]
 \begin{picture}(200,70)(25,6)
 \put(140,30){\circle{12}}
 \put(152,27){$2g-2$}
 \put(140,36){\line(0,1){20}}
 \put(58,56){\circle*{7.5}} \put(55,66){$b$}
\put(40,66){$\recip{r}\Bigl($}
 \put(59,56){\line(1,0){21}}
 \put(85,56){\circle*{7.5}} \put(80,66){$2b$}
 \put(94,56){\dots}
 \put(115,56){\circle*{7.5}} \put(101,66){\dots}
 \put(117,56){\line(1,0){21}}
 \put(140,56){\circle*{7.5}} \put(122,66){$(r-b)b$}
 \put(141,56){\line(1,0){21}}
 \put(165,56){\circle*{7.5}}
 \put(174,56){\dots}
 \put(195,56){\circle*{7.5}} \put(180,66){$r-b$}
\put(206,66){$\Bigr)$}
 \put(114,30){\line(1,0){20}}
 \put(113,30){\circle*{7.5}}
 \put(122,10){\line(1,1){15}}
 \put(121,8){\circle*{7.5}}
 \put(89,8){\dots}
 \put(107,8){\line(1,0){28}}
 \put(141,8){\dots}
 \end{picture}
 \caption{The quasistellar graph $\Ga(g,\sB)$.}
 \label{fig!star}
 \end{figure}
The central vertex is a divisor $B$ of self-intersection $B^2=2g-2$
(ineffective if $g=-1$, that is, $B^2=-4$); for each term
$\recip{r}(a,r-a)\in\sB$, it meets the $b$th curve in a chain of $r-1$
exceptional $-2$-curves, where $ab\equiv1$ mod $r$. (Compare Belcastro
\cite{Be} for the lattice of the famous 95 families of K3
hypersurfaces.)

The polarising $\Q$-divisor $D=f^*D=B+\sum m_iE_i$ has exceptional
curves in each chain weighted by the monotone sequence of arithmetic
progressions
 \[
 \recip{r}\Bigr(
 b,2b,\dots,b(r-b-1),b(r-b),(b-1)(r-b),\dots,2(r-b),r-b \Bigl),
 \]
giving
 \[
 D\cdot E_i=0, \quad\hbox{and}\quad
 D^2=D\cdot B=2g-2+\sum\frac{b(r-b)}{r}\,.
 \]

The lattice $L=L(g,\sB)$ has signature $(+1,-\sum(r-1))$ and
discriminant $(\prod_\sB r)D^2$. Although $L$ is determined up to
isomorphism by the numerical data $g$ and $\sB$, its embedding into the
standard K3 lattice
 \[
 L\into H^2(T,\Z)=\LK
 \]
is not in general completely determined by $g$ and $\sB$, and is an
additional topological invariant of $S,D$. For example, if $S$ is
nonsingular then $L=\Z\cdot D=\Span{2g-2}$ is the lattice with one
generator of square length $2g-2$, but it can happen that $D$ is
divisible in $\Pic S$ by some integer $n$ with $n^2\divides g-1$, and
then $L\subset\LK$ is not a primitive sublattice. If $L$ has high rank,
it can have several inequivalent primitive embeddings $L\into\LK$,
having nonisomorphic orthogonal complements.

The Hilbert series of a Fano 3-fold $V$ is controlled by the same
combinatorics: the lattice $L(g,\sB)$ is still there (even if $g=-2$).
For simplicity, assume that $V$ has only terminal quotient
singularities $\recip{r}(1,a,r-a)$. Make the economic resolution $Y\to
V$ of these singularities by successive $(1,a,r-a)$ weighted blowups.
These are the Kawamata blowups of Corti, Pukhlikov and Reid
\cite{CPR}, 3.4.2 (see also \cite{YPG}, 5.7; but note that the
chain of Kawamata blowups is a unique ``best'' choice of economic
resolution that came on line after \cite{YPG}). The lattice $L(g,\sB)$
is the Picard lattice $\Pic Y$, with bilinear product
 \[
 (D_1,D_2)\mapsto (-K_Y)\cdot D_1\cdot D_2.
 \]
 If $g\ge-1$ and we assume that $V$ (or its deformation) has a K3
elephant, then $L(g,\sB)$ is a sublattice of $\LK$. If $g=-2$, so that
$B^2=-6$, we have to imagine that we are looking at a nonexistent K3
surface having a polarising divisor $D$ with $h^0(D)=-1$. In this case,
there is no obvious a priori relation between $L(g,\sB)$ and $\LK$;
however, as we said in \ref{sssec!ineq}, there are probably only a few
cases, and it is not a substantial restriction for a lattice in this
range to be a sublattice of $\LK$, so who knows?

 \section{From K3s to Fano 3-folds} \label{sec!proj}
 Just as the general theory of~\ref{ssec!ele}, what we said in
Section~\ref{sec!Hi} applies to K3 surfaces and to Fano 3-folds alike.
We said in \ref{sssec!ineq} that K3 surfaces go on for ever and for
ever, with many baskets continuing to exist for infinitely many values
of $g$. However, as in the more specific 3-fold theory of
\ref{ssec:specific}, there are results of Kawamata
\cite{Ka1}--\cite{Ka2} saying that Fano 3-folds are bounded. By
analogy with Theorem~\ref{th!clF}, we guess that only a couple of
thousand families exist, and we eventually aspire to a precise
classification. Most Fano \hbox{3-folds} have $H^0(-K_V)\ne0$, and, as
a first attempt, we can take those having a K3 elephant as typical,
and allow ourselves extra generality assumptions on $V$ such as
$\Q$-smoothing, Brill--Noether general behaviour of linear systems
(nonhyperelliptic, etc.).

 \subsection{Takagi's results} \label{sssec!F3}
 In his lecture at the Singapore conference \cite{T1}, Takagi sketched
a preliminary classification of prime $\Q$-Fano \hbox{3-folds} of genus
$g\ge2$ having at worst terminal quotient singularities and a K3
elephant. The assumption that $|{-}K_V|$ contains a K3 is a strong
condition: it guarantees that there is a resolution on which $-K_{V'}$
is nef. In this case, the anticanonical system $|{-}K_{V'}|$ is a big
linear system of K3 surfaces, and defines a birational model of $V$ as
an indecomposable Gorenstein Fano 3-fold $X$ (Fano would appreciate
this picture at once!). Note that the anticanonical system $|{-}K_V|$
itself (but not $|{-}2K_V|$, $|{-}3K_V|$, etc.)\ is invariant under
passing to the economic resolution of points of index $r\ge2$. Mukai's
Theorem~\ref{th!clF} applies to $X$, thus reducing the study of
\hbox{$\Q$-Fano} \hbox{3-folds} to problems in projective geo\-metry
concerning sections of homogeneous spaces with special configurations
of planes and singularities. In particular, Takagi proves that if $V$
has at least one terminal quotient singularity then $g\le8$; the cases
$g=7$ and $g=8$ is accessible since $X=\fie_{-K_V}(V)$ is a linear
section of a Grassmann or orthogonal Grassmann variety. Takagi's thesis
\cite{T} also settled all cases with only $\half(1,1,1)$ singularities.

 \subsection{Low codimension} \label{sssec!lowc}
 We study Fano 3-folds in terms of the (pluri-) anticanonical ring
 \[
 R=R(V,-K_V)=\C[x_0,\dots,x_N]/I_V,
 \]
 where $x_0,\dots,x_N$ (more usually $x,y,\dots,u$, etc.)\ are
homogeneous generators of weights $a_0,\dots,a_N$ and $I_V$ is the
ideal of all relations. Since $-K_V$ is ample, $\Proj R$ is the
anticanonical model of $V$ in w.p.s.:
 \[
 V=\Proj R\subset\PP^N(a_0,a_1,\dots,a_N).
 \]
 Here we take the complete anticanonical ring $R(V,-K_V)$, so that
$V\subset\PP^N$ is projectively normal; it is known that $V$ is
projectively Gorenstein, that is, $R$ is a Gorenstein ring (see for
example Goto and Watanabe \cite{GW}).

The {\em codimension} of $R$ means the codimension $N-3$ of the
anticanonical ring $R(V,-K_V)$. It is a measure of the difficulty of
studying $V$ directly. In codimension~2, a Gorenstein ring is a c.i.\
(Serre), and in codimension~3, it is given by the $2k\times2k$ diagonal
Pfaffians of a skew $(2k+1)\times(2k+1)$ matrix (Buchsbaum and Eisenbud);
only $5\times5$ occurs here, with the single exception of the classic
case $Q_1\cap Q_2\cap Q_3\subset\PP^6$.

The cases when $V$ is a hypersurface or codimension~2 c.i.\ are settled
in Fletcher \cite{Fl}, and the codimension~3 cases in Alt{\i}nok
\cite{A}. There are 95 families of Fano 3-fold hypersurfaces, 84
codimension~2 c.i.s (including Fletcher's example mentioned in
\ref{ssec!uh} with $H^0(-K_V)=0$), and 70 families of codimension~3
Pfaffians. These varieties are given by explicit equations, that we
can write down without reference to K3s, and depend on an irreducible
parameter space. Their construction is unobstructed, and every family
of K3 surfaces extends to a family of Fano 3-folds.

 \subsection{Weighted Grassmannian} \label{ssec!wGr}
 In codimension~3, the five Pfaffians are a weighted form of the equations
of the Grassmannian $\Grass(2,5)$ (see \ref{rmk!Pf}). Thus a
codimension~3 K3 surface $S\subset\PP(a_0,\dots,a_5)$ has Hilbert
series
 \[
 1-\sum t^{b_i} + \sum t^{k-b_j} - t^k,
 \]
where $k=\sum_0^5 a_i$ and $2k=\sum_1^5 b_i$, so that the skew matrix
in \ref{rmk!Pf} has entries of weight $\wt x_{ij}=k-b_i-b_j$ (compare
Reid \cite{Ki}, Example~3.8). We can view the $x_{ij}$ as coordinates
on $\bigwedge^2U\otimes L$, where $U$ is a weighted $\C^5$ with weights
$c-b_i$, and $L=\C$ with weight $k-2c$ (for some convenient $c$); then
the equations of $S$ are the pullback of the weighted Grassmann
$\wGr\subset\PP(\bigwedge^2U)$ (tensors of rank~2) by a weighted
projective map $\PP(a_0,\dots,a_5)\broken \PP(\bigwedge^2U)$. In the
spirit of Mukai's work discussed in \ref{ssec:specific}, the interesting
question is to construct directly an exceptional (orbi-) vector bundle
$\sE$ on $S$ whose Serre module $H^0_*(\sE)$ is generated by 5 sections
$u_i\in\sE(c-b_i)$, giving the embedding
$S\into\wGr\subset\PP(\bigwedge^2U)$.

The point here is to get away from the algebraic treatment of the
equations of $S$ to a more geometric understanding of what they mean;
that is, eventually, to replace the appeal to Buchsbaum and Eisenbud's
algebraic result on codimension~3 Gorenstein rings by arguments in
the style of Mukai in terms of exceptional vector bundles on the K3
section. There are preliminary results in this direction due to Corti
and Reid \cite{CR}. This raises important unsolved problems in higher
codimension: weighted versions of Mukai's symmetric spaces make sense,
but seem to give few examples of varieties of small coindex.

 \subsection{Gorenstein unprojection}
 In codimension~$\ge4$, there is no analogous structure theorem for
Gorenstein rings; our current strategy for K3s and Fano \hbox{3-folds}
in codimension~4, 5, etc., is based (in most cases) on {\em Gorenstein
projection}. A model case is the projection $S_d\broken S_{d-1}$ of a
del Pezzo surface $S_d\subset\PP^d$ from $P\in S_d$ (for simplicity,
think of $d\ge4$). Here both $S_d$ and $S_{d-1}$ are anticanonical: the
blowup $S'\to S_d$ of $P$ has discrepancy~1, so that the anticanonical
map $S'\to S_{d-1}$ given by $|{-}K_{S_{d-1}}|=|{-}K_{S_d}-E|$ is the
same thing as the linear projection from $P$. 

Although the study of a Fano 3-fold $V\subset\PP^N(a_0,a_1,\dots,a_N)$
is a biregular problem in the first instance, we can use birational
methods to attack it. Fano's linear projection of
$V=V_{2g-2}\subset\PP^{g+1}$ with $g\ge7$ from a line can also be
interpreted as making a blowup $V_1\to V$, then recalculating the
anticanonical ring $R(V_1,-K_{V_1})$ as a subring of $R(V,-K_V)$, thus
deducing a birational map $V\broken V'=\Proj R(V_1,-K_{V_1})$. In
Fano's study, $V$ is the unknown, and is hard to work with because it
has high codimension, whereas its projection $V'$ has smaller
codimension, so it, together with the exceptional divisor of the
projection $F\subset V'$, may be more tractable.

For us, the key point is that the birational relation between $V$ and
$V'$ (geo\-metrically a projection) can be handled in terms of
inclusions
 \[
 R(V',-K_{V'})\subset R(V,-K_V)
 \]
 between Gorenstein rings. This area has recently been clarified by
Papadakis and Reid's treatment of the inverse birational map
$V'\broken V$ as {\em Gorenstein unprojection} or ``constructing
big Gorenstein rings from small ones'' (see \cite{PR} and \cite{Ki},
and compare Kustin and Miller \cite{KM}); in \cite{PR}, the inverse is
constructed in terms of the adjunction formula for the
Grothendieck--Serre dualising sheaf.

 \subsection{Type~I unprojection} \label{ssec:TI}
 We discuss here the most straightforward case of projection,
corresponding to Kustin and Miller unprojection, which already applies
to the majority of K3s in codimension $\ge4$. Let
 \[
 S\subset\PP(a_1,\dots,a_N) \quad\text{or}\quad
 V\subset\PP(1,a_1,\dots,a_N)
 \]
 be a K3 surface or Fano 3-fold. We only treat the Fano case here,
leaving the reader to make the obvious modifications in the K3 case.

A coordinate point $P_k=(0,\dots,1,\dots,0)\in V$ is a {\em Type~I
centre} if $P_k\in V$ is a terminal quotient singularity
$\recip{r}(1,a,r-a)$, where $\wt x_k=r$, and global coordinates of the
projective space $x_0,x_i,x_j$ of weight equal to $1,a,r-a$ provide
local orbifold coordinates at $P_k$. Let $V_1\to V$ be the Kawamata
blowup, that is, the weighted blowup of $P_k$ with weights $1,a,r-a$;
then the exceptional locus $E\iso\PP(1,a,r-a)$ has minimal discrepancy
$\recip{r}$, and, following \cite{CPR}, we write
 \[
 A=-K_V \quad\hbox{and}\quad B=-K_{V_1}=A-\recip{r}E.
 \]
 Now because $\wt(x_0,x_i,x_j)=(1,a,r-a)$, and $x_0,x_i,x_j$ vanish
along $E$ with multi\-plicity $\recip{r},\frac{a}{r},\frac{r-a}{r}$
(see \cite{CPR}, Proposition~3.4.6), they belong to the subring
$R(V_1,B)$. Thus $R(V_1,B)$ is the subring
$k[x_0,\dots,\widehat{x_k},\dots,x_N]\subset R(V,A)$ obtained by
eliminating $x_k$ only. Then
 \[
 V'=\Proj R(V_1,B)\subset\PP(1,a_1,\dots,\widehat{a_k},\dots,a_N)
 \]
 is a weak Fano 3-fold containing the plane $\Pi=\PP(1,a,r-a)$. Be
warned that $V'$ is not a Mori Fano 3-fold, since the Weil divisor
$\Pi$ is not in $\Z\cdot(-K_{V'})$; one usually expects $V'$ to be the
{\em midpoint} of a Sarkisov link, compare \cite{CPR}, 4.1, (3).

 A simple example of a Type~I unprojection was treated in
Example~\ref{exa!2/3}. Type~I projections include the construction of
the 64 quadratic involutions of \cite{CPR}, 4.4 (or rather, of the
first half $X\broken Z$, up to the midpoint of the link).

 \begin{exc} \label{ex!t1}
 In the numerical data $g,\sB$ for a K3 surface
(Section~\ref{sec!Hi}), replace an element $\recip{r}(a,r-a)$ of $\sB$
by two elements $\recip{a}(r,-r)$ and $\recip{r-a}(r,-r)$, and assume
that $D^2$ remains positive. Study how this numerical projection
affects the rhs of (\ref{eq!K3-s}): show that it subtracts
$\frac{t^r}{(1-t^r)(1-t^a)(1-t^{r-a})}$ from the Hilbert series
$P_{S,D}(t)$, and reduces $D^2$ by $\recip{ra(r-a)}$. Compare
\cite{PR}, Exercise~2.7, and see Example~\ref{exa!dbeg1},
page~\pageref{Altinok3(63)} for a numerical instance.
 \end{exc}

 The effect of a Type~I projection as in Exercise~\ref{ex!t1} on the
numerics of Section~\ref{sec!Hi} gives the set of numerical data of K3
surfaces the structure of a directed graph, with comparatively few
connected components. With few exceptions, families of K3 surfaces in
codimension~$\ge2$ have projections to smaller codimension, most
commonly of Type~I. Of the 142 codimension~4 K3s in the K3 database,
116 have a numerical Type~I centre. All but 2 of the remaining cases
are covered by the higher types of projection discussed in Reid
\cite{Ki}, Section~9 and \cite{T4}. These higher projections, and the
small core of exceptions not admitting any projections, are interesting
and demand further study; we suspect that these more complicated K3s
are unlikely to extend to Fano 3-folds.

 \subsection{Tom and Jerry unprojections to codimension 4} \label{T&J}
 When applicable, a Type~I projection as described in \ref{ssec:TI}
reduces the study of a Fano $V$ in w.p.s.\ $\PP^N$ to a variety $V'$ in
w.p.s.\ $\PP^{N-1}$, but specialised to contain a weighted projective
subspace:
 \begin{equation}
 \PP(1,a,r-a)\subset V'\subset \PP^{N-1}.
 \label{eq!t&j}
 \end{equation}
 Thus Fano \hbox{3-folds} whose numerical data admits a codimension~4
candidate can often be studied via projections. However, setting up the
unprojection data (\ref{eq!t&j}) is still a difficult problem;
complicated features of $V$, such as obstructed equations or reducible
moduli spaces, must be faithfully reproduced in the unprojection data.

In (\ref{eq!t&j}), the w.p.s.\ $\PP(1,a,r-a)$ is a codimension~4 c.i., and
$V'$ is a codimension~3 variety given by the Pfaffians of a $5\times5$
skew matrix. Thus to obtain codimension~4 Fanos with this type of
projection, the problem is how to put a codimension~4 c.i.\ inside a
$5\times5$ Pfaffian. Similarly for K3s.

There are two different solutions to this problem, called {\em Tom} and
{\em Jerry}, treated in detail in Papadakis \cite{P}--\cite{P1}
(see also
\cite{Ki}, Examples~6.4 and 6.8, and Section~8). In these two cases, we
specialise the skew $5\times5$ matrix defining $V'$ to
 \[
 M_{\mathrm{Tom}}=
 \left(
 \begin{array}{c@{\enspace}c@{\enspace}ccc}
 x_{12} && x_{13} & x_{14} & x_{15} \\[2pt]
 \cline{2-5}
 &\vline&a_{23}&a_{24}&a_{25} \\
 &\vline&&a_{34}&a_{35} \\
 &\vline&&&a_{45}
 \end{array}
 \right)
 \quad\hbox{or}\quad
 M_{\mathrm{Jerry}}=
 \left(
 \begin{array}{cccc}
 a_{12}&a_{13} & a_{14} & a_{15} \\
 &a_{23}&a_{24}&a_{25} \\[2pt]
 \cline{2-4}
 &&x_{34}&x_{35} \\
 &&&x_{45}
 \end{array}
 \right)
 \]
where the 6 entries $a_{ij}$ in the bottom left $4\times4$ block of
$M_{\mathrm{Tom}}$ specialise to lie in a codimension~4 complete
intersection ideal $(y_1,\dots,y_4)$ (and ditto for the 7 entries
$a_{ij}$ in the first two rows and columns of $M_{\mathrm{Jerry}}$).
In either case, the theoretical construction of \cite{PR} or \cite{KM}
gives the unprojection $V\subset\PP^N$ and its anticanonical ring, and
Papadakis \cite{P1}, Section~5 gives an explicit presentation of the
ring.

Tom and Jerry occur in hundreds of constructions of Gorenstein
codimension~4 rings with $9\times16$ presentation, and seem to be
related to the respective cones over $\PP^2\times\PP^2$ and over
$\PP^1\times\PP^1\times\PP^1$ and their weighted homogeneous
deformations; but exactly what this means remains to be elucidated. See
\cite{Ki}, Section~8 for a more detailed discussion.

To deal with our codimension~4 K3s and Fanos, we have to determine
whether Tom and Jerry matrixes can be set up to give $V$ having only
Mori category singularities. There are a few hundred problems here, and
we have only just started working systematically on the Fanos. For
each case where the numerical data of the K3 admits a Type~I
projection, we can ask for a Tom or Jerry matrix -- from the
experience of Alt{\i}nok's thesis \cite{A} and \cite{CPR}, we suspect
that in each case at least one of Tom or Jerry exists, and that when
both exist, the two families give isomorphic polarised K3s, so do not
correspond to different irreducible components of moduli.

When we ask the same question for Fanos 3-folds, one of several things
may happen: either of the Tom and Jerry matrixes could exist to give a
Mori Fano \hbox{3-fold}, or either could fail. We suspect that each of
these possibilities happens in many numerical cases. When both Tom and
Jerry exist, the resulting unprojected varieties may give essentially
different different Fano 3-folds: see \cite{Ki}, Examples~6.4 and~6.8,
based on Takagi's thesis, where Tom and Jerry unprojections with the
same numerical data correspond to Fano 3-folds that are not biregular,
and not obviously birational. Since we have several hundred cases to
settle, and much of the calculation comes down to checking that
monomials of suitable degree exist, we hope that most of the
calculation can eventually be entrusted to the computer.

 \section{The K3 database in Magma}
 The program outlined in Sections~\ref{sec!Hi}--\ref{sec!proj}
envisaged listing the many thousand possible values $g,\sB$ for the
numerical data of families of K3s or Fano \hbox{3-folds}. For each
$g,\sB$ in the list, we calculate the Hilbert series and try to deduce
a plausible shape for the anticanonical ring and its possible
pro\-jections. In some cases, we can establish the existence of a
quasi-smooth surface, extension to Fano 3-fold, and connected
components of moduli -- these are all questions that become nontrivial
when $R(S,D)$ has codimension $\ge4$. These calculations are akin to
traditional work on embedding curves and surfaces in Italian algebraic
geometry, although rational contributions from quotient singularities
add a certain spice, and unlimited possibilities for errors of
arithmetic. With patience, any particular calculation can be done by
hand, but automation has obvious advantages when doing these
calculations on an industrial scale.

 Several items of this program have already been carried out for
polarised K3 surfaces. Hilbert series methods apply to many other
situations in algebraic geo\-metry involving graded rings or modules;
current work in progress by students at Warwick includes Suzuki \cite{S}
on Fano 3-folds of Fano index $f\ge2$, Keenan's project \cite{Ke} on
subcanonical curves and Buckley's study \cite{Bu} of polarised
Calabi--Yau \hbox{3-folds} with strictly canonical singularities.
Nonetheless, rings over K3 surfaces have several advantages that make
them an ideal target for a computer study: the numerical data and the
RR statement are simple to state, there is no cohomology, and most K3
surfaces have projections to smaller codimension. Moreover, work over
the last 20 years (Fletcher \cite{Fl}, Alt{\i}nok \cite{A} and Corti,
Pukhlikov and Reid \cite{CPR}) already provides us with several
hundred worked out examples of what to expect.

 We have programmed these calculations as routines in John Cannon's
computer algebra system Magma \cite{Ma}. The results of all the
calculations in codimension $\le4$ are assembled as a database
containing 391 families of K3 surfaces. We describe here a little of
what our code does and how to use it. Although more sophisticated
upgrades are projected (and implemented in part, see the Graded Ring
Database Website \cite{Br1}), we restrict ourselves to the first
working version of our routines, in the widely available export~2.8 of
Magma \cite{Ma}. Please see \cite{Br1} for more recent work; this
includes a prototype database of K3s up to codimension~10, and will
document upgrades as they come on line.

 \begin{exa}[making a single surface]\label{exa!dbeg1}
 Suppose we want to study a polarised K3 surface $S,D$ with
ineffective $D$ (that is, $g=-1$, see \ref{ssec!uh}), and having basket
of singularities $\sB=\{\recip{2}(1,1),\recip{5}(1,4),
\recip{13}(3,10)\}$. As usual, {\tt >} at the beginning of a line is
the Magma prompt. Anything after this prompt is user input, while
anything else is Magma output (occasionally subjected to a little
editing for legibility). We have formatted the Magma input and output
along the lines of our mathematical writings, and you should not have
too much trouble interpreting the language.
\begin{verbatim}
> Q := Rationals();
> R<t>:=PolynomialRing(Q); // Omit these lines at your peril.
> B := [ [2,1], [5,1], [13,3] ];
> S := K3Surface(-1,B); // We input genus -1 and Basket B.
> S;
Codimension 4 K3 surface with data
  Weights: [ 3, 4, 5, 6, 7, 10, 13 ]
  Numerator: t^48 - t^36 - t^35 - t^34 - t^33 - t^32 - t^31 - t^30
    + t^27 + 2*t^26 + 2*t^25 + 2*t^24 + 2*t^23 + 2*t^22 + t^21
    - t^18 - t^17 - t^16 - t^15 -  t^14 - t^13 - t^12 + 1
  Basket: [ 2, 1 ], [ 5, 1 ], [ 13, 3 ]
\end{verbatim}
The computer has calculated the Hilbert series, and, after
experimenting with a number of possibilities, has found
 \[
 (1-t^3)(1-t^4)(1-t^5)(1-t^6)(1-t^7)(1-t^{10})(1-t^{13})
 \]
 as a plausible denominator, corresponding to generators for the ring
$R(S,D)$, or an embedding of $S$ in w.p.s.\
$\PP=\PP^6(3,4,5,6,7,10,13)$. The printout lists the conclusion: a
candidate family of codimension~4 K3 surfaces $S\subset\PP$ with the
stated Hilbert numerator. We should be clear that in general
 \begin{quote}
{\em there is no a priori guarantee that this or any other candidate
surface proposed by our Magma functions exists as a subvariety in the
indicated w.p.s.\ with the indicated properties.}
 \end{quote}
After all, at this stage the computer has done nothing more than a
formal game with denominators for Hilbert series. Having said that,
assume for a moment that $S\subset\PP^6(3,4,5,6,7,10,13)$ really exists,
polarised by the divisorial sheaf $\Oh_S(D)=\Oh_S(1)$. Clearly $D$ is
not an effective divisor, since no $x_i$ has weight $a_i=1$, but its
multiples do exist: in particular, $S$ has a single effective divisor
in $|3D|$, $|4D|$, $|5D|$, a pencil $|6D|$, etc.

The {\tt Numerator} is a polynomial $p=P(t)\prod(1 - t^{a_i})$, where
$P(t)$ is the Hilbert series of $R(S,D)$
 \[
 P_S(t)=\sum_{n\ge0} h^0(S,nD)t^n = \frac{p}{\prod(1 - t^{a_i})}\,.
 \]
Here the product is taken over the weights $a_i$ of the w.p.s. At first
sight, the 7 negative terms $-t^{12}-t^{13}-\cdots$ in the numerator
output by Magma might suggest that $S\subset\PP$ could have 7 equations
in degrees 12, 13, 14, 15, 16, 17, 18. Gorenstein rings with
$7\times12$ resolution certainly exist, so this is not completely
excluded, but this structure appears very rarely (compare Papadakis and
Reid \cite{PR}, 2.8); perhaps this candidate appears more exotic than it
actually is. The $9\times 16$ format is ubiquitous, and experience
suggests that this is really a $9\times 16$ resolution in disguise.

How can we prove that this surface actually exists? The method of
Alt{\i}nok's thesis \cite{A} is to look for a codimension~3 K3
surface that is already known to exist, and could be a projection of
the desired surface. Then we could go on to show that the
codimension~3 surface can be made to contain an {\em unprojection
divisor\/} as in \ref{ssec:TI} that can be contracted or {\em
unprojected\/}, to give the codimension~4 surface. This second stage
is the hard part, and we pass over it for the present (see \cite{PR},
\cite{P1} and \cite{Ki}), although in time the computer will have
something useful to add here too. One point is that as we write there
are many kinds of projection for which the unprojection is still
something of a mystery (see
\cite{Ki}, Section~9). That is, the {\em theory} of unprojection is not
quite complete enough yet to handle all the cases we might need.

However, the K3 database is very good at finding projections of the
types we already understand, such as the Type~I projections of
Section~\ref{sec!proj}. These occur when we project from a cyclic
quotient singularity $\recip{r}(a,r-a)$ whose local orbifold
coordinates come from generators of $R(S,D)$ of weight $a,r-a$; a
singularity with this property is a {\em Type~I centre}. In our
example, we can see a possible Type~I centre, the singularity
$\recip{13}(3,10)$. The point is that the local weights 3 and 10, a
priori only defined mod~13, can arise as the global generators
corresponding to the weights 3 and 10 of the w.p.s. All the other
singularities fail this test, since the w.p.s.\ has no weight $1$.
\begin{verbatim}
> DB := K3Database("t"); /* Loads the DB, requiring "t" as
 variable in Hilbert series. This takes several seconds. */
> Centres(~S,DB); // Searches DB for projection centres of S.
> S;
Codimension 4 K3 surface with data
  Weights: [ 3, 4, 5, 6, 7, 10, 13 ]
  Numerator: t^48 - t^36 - t^35 - t^34 - t^33 - t^32 - t^31 - t^30
 + t^27 + 2*t^26 + 2*t^25 + 2*t^24 + 2*t^23 + 2*t^22 + t^21
 - t^18 - t^17 - t^16 - t^15 -  t^14 - t^13 - t^12 + 1
  Basket: [ 2, 1 ], [ 5, 1 ], [ 13, 3 ]
  Centre 1: [ 5, 1, 4 ] has Type 2 projection to 10 in codim 1
  Centre 2: [ 13, 3, 10 ] has Type 1 projection to 42 in codim 3
\end{verbatim}
This says that the singularity $\recip{13}(3,10)$ is a Type~I centre
as expected. It has also found a Type~II centre, that we pass over for
the moment. But it does more: Magma applies the calculus of
Exercise~\ref{ex!t1} to predict the weights and basket of the image of
a Type~I projection, and searches the database for surfaces with the
right properties, finding the K3 surface numbered 42 in the database
as a plausible image of the projection from this centre. (The internal
numbering of items in the database is arbitrary, and may differ from
session to session.) We ask what it is:
\label{Altinok3(63)}
\begin{verbatim}
> K3Surface(DB,42);
Codimension 3 K3 surface, number 42, Altinok3(63), with data
  Weights: [ 3, 4, 5, 6, 7, 10 ]
  Numerator: -t^35 + t^23 + t^22 + t^21 + t^20 + t^19
 - t^16 - t^15 - t^14 - t^13 - t^12 + 1
  Basket: [ 2, 1 ], [ 3, 1 ], [ 5, 1 ], [ 10, 3 ]
\end{verbatim}
Now if we know that the image codimension~3 surface
$S'\subset\PP(3,4,5,6,7,10)$ exists, and contains the stratum
$\PP(3,10)\subset\PP$ as unprojection divisor, then it has the
required geometric properties to be unprojected, and we conclude from
\cite{PR} that the original surface exists. Indeed, since we know a
lot about unprojection, we can even write down the equations of the
codimension~4 surface $S$ based on those of the codimension~3 surface
$S'$ and its unprojection divisor $\PP(3,10)$. (Since $\PP(3,10)$ is
the c.i.\ defined by $x_4=x_5=x_6=x_7=0$, and the unprojection
variable $s$ has weight $13$, we verify that $S$ also needs equations
in degree $19,20$, so our intuition that $S$ has a $9\times16$
resolution was right.)

The basic existence part is easy. The structure theorem for Gorenstein
rings of codimension $\le3$ makes it easy for us to figure out
equations. Alt{\i}nok goes further to check that a surface containing
the unprojection divisor can be made. Using Alt{\i}nok's numbering in
\cite{A1}, Tables~5.1 and~5.2, the codimension~3 surface is
\Alt{3}{63}, as the database tells us, while the codimension~4 surface
is \Alt{4}{104}.

Of course, the projection calculus is easy to do by hand: we could
simply have searched the lists of \cite{Fl} and \cite{A} to find this
surface (as we always did in the dark ages before the great Graded Ring
Database Revolution). For that matter, we could have just called for a
surface having the desired basket and again found this one without
trouble. At this stage, we have not gained much by having the results
of these calculations in a database. The next section contains a
calculation for which the database is essential.
 \end{exa}

 \begin{exa}[searching the database]
 As databases go, our current database of 391 K3s is tiny. Even so,
its key advantage over a printed list is that we can search it quickly
and accurately. Suppose, for example, that we are interested in finding
out whether our lists contain a K3 having a singularity of index~17.
In this context, 17 is a large number, taking the sublattice
$L(g,\sB)\subset\Pic S$ close to the maximum that can live in the K3
lattice $\LK$ (see \ref{sssec!K3L}). To do this search in Magma, we
first summon the database as usual:
\begin{verbatim}
> DB := K3Database("t");
> #DB; // # returns the number of elements of a list.
391 \end{verbatim}
Singularities are denoted {\tt p = [r,a]}, with $a$ coprime to the
index {\tt r = p[1]} and $a<r/2$. So we must search the database for
surfaces having a singularity $[r,a]\in\sB$ having $r=17$.
\begin{verbatim}
> surfaces := [ S : S in DB | &or[ p[1] eq 17 : p in Basket(S) ] ];
> #surfaces;
2
\end{verbatim}
Note that {\tt \&or} taken over a list of Boolean values is equivalent
to ``there exists $p$ in $\sB$ with $p[1]$ equals 17''. The last line
says that there are just two surfaces in DB with a singularity of
index~17. We ask to see them.
\begin{verbatim}
> for S in surfaces do print S;
> end for;
Codimension 2 K3 surface, number 1, Fletcher2(82), with data
  Weights: [ 3, 4, 7, 10, 17 ]
  Numerator: t^41 - t^21 - t^20 + 1
  Basket: [ 2, 1 ], [ 17, 7 ]
Codimension 4 K3 surface, number 250, Altinok4(79), with data
  Weights: [ 2, 3, 5, 5, 7, 12, 17 ]
  Numerator: t^51 - t^41 - t^39 - t^37 - t^36
   + t^29 + t^27 + t^26 + t^25 + t^24 + t^22
   - t^15 - t^14 - t^12 - t^10 + 1
  Basket: [ 17, 5 ]
\end{verbatim}

Again, we want to know that these surfaces exist. The codimension~2
surface is Fletcher \cite{Fl}, List~13.8, no.~82, and it is easy to
write down its equations.

For the codimension~4 surface, the technique of the previous example
locates an image of projection from the Type~I centre
$\recip{17}(5,12)$ in the database, and we can construct $S$ by
unprojection as before. This codimension~4 example $S$, called
\Alt{4}{79} is interesting in that it has exactly the same Hilbert
series as the codimension~4 c.i.\ $T_{10,12,14,15}$; however, none of
the equations of the c.i.\ $T$ can involve the last variable $u_{17}$,
so that $T$ is a kind of weighted cone over a curve. Deformations of
polarised K3 surfaces are unobstructed, as are c.i.s, so that having
proved that $S$ exists, we have found two different irreducible
components of the Hilbert scheme, one containing $S$, the other
containing $T$.
 \end{exa}

\subsection{K3 database functions in Magma}

As with much computer algebra, just a few internal functions do the
work, and a lot of the rest is cosmetic renaming to simplify the
user's access to these. We give a brief and possibly inadequate
description of some core functions to give some idea of what is
happening inside the computer; these remarks can also serve as a first
introduction to the code for anybody wishing to modify it for use in
other contexts.

 \subsubsection{Listing all baskets} The function
\begin{verbatim}
> Baskets(n);
\end{verbatim}
returns a list of all possible baskets $\sB=\{\recip{r}(a,r-a)\}$ for
K3 surfaces of singular rank $\sum(r-1)<n$. For example, the following
commands generates a list of all baskets with $\sum(r-1) = 11$, and
asks for the 44th basket in the list.
\begin{verbatim}
> BB, gg := Baskets(12); /* BB is the list of baskets,
   gg the parallel list of minimal genera */
Actual number: 329 ... Checking degrees ...
> BB11 := [ B : B in BB | not (B eq [])
    and &+[p[1]-1 : p in B] eq 11 ];  // &+ is sum over set.
> #BB11;
109
> b := BB11[44];
> b, gg[Index(BB,b)];
[ [ 4, 1 ], [ 5, 2 ], [ 5, 2 ] ]
0
> Degree(-1,b);   
-17/20
> Degree(0,b); 
23/20
\end{verbatim}

Recall from \ref{sssec!ineq} that the possible baskets $\sB$ on a K3
surface are limited by two inequalities on the pairs $r,a$. The first
is given by the rank of the Picard group: $\sum(r-1)<n=20$. The
argument $n$ of {\tt Baskets(n)} is the bound (we eventually set
$n=20$ for complete lists). The other inequality, saying that the
polarised surface $S,D$ has degree $D^2>0$, involves the genus $g$ of
$S$. If $g=-1$ or $0$, a basket may give $D^2\le0$. The second return
value of {\tt Baskets} (called $gg$ in the first line above) is a
parallel sequence of minimal genera that give $D^2>0$. In the above
example, we checked that $g=-1$ is illegal for the 44th basket,
because it gives $D^2=-17/20$.

The numerical data $g,\sB$ with $D^2>0$ determines a Hilbert series
$P(t)=\sum P_nt^n$. If it actually corresponds to a K3 surface,
then vanishing implies that $P_n\ge0$ for all $n\ge1$. However, this
does not follow from the inequality $D^2>0$. Indeed, there exist
numerical data $g,\sB$ with a negative coefficient $P_n$; curiously,
there are just three of these eccentrics (compare \cite{CPR}, 7.9).
The function {\tt Baskets} also checks that $P_n\ge0$ for the first 20
coefficients before confirming the minimum genus $g$. When $n=20$, the
result is all 6640 possible baskets, together with, for each, the
minimum genus that makes the Hilbert series positive.

 \subsubsection{{\tt HilbertSeries(g,B)}, etc.} The Magma function
\begin{verbatim}
> HilbertSeries(g,B);
\end{verbatim}
has the formula (\ref{eq!K3-s}) built in, and simply evaluates it at
the data $g,\sB$ as a rational function in $t$. For example,
\begin{verbatim}
> R<t>:=RationalFunctionField(Q); // Require "t" as variable.
> P:=HilbertSeries(-1,[[ 3, 1 ], [ 4, 1 ], [ 11, 2 ]]);
> P*(1-t^2)*(1-t^3)*(1-t^4)*(1-t^11);        
-t^20 - t^15 - t^13 - t^11 + t^9 + t^7 + t^5 + 1
\end{verbatim}
suggests the candidate surface $S\subset\PP(2,3,4,5,7,9,11)$ with the
stated $g$ and $\sB$. (To make this work for other types of graded
rings, you need to figure out what Hilbert series you want to use, and
program it in as a substitute for our formula (\ref{eq!K3-s}).)

In the above example, we put in the denominator by hand, and this
turned out to be a lucky guess. The heart of the whole package is a
suite of functions to make a reasonable analysis of this Hilbert
series, using a few tricks based on the experience of Alt{\i}nok's
thesis \cite{A}. One expects generators corresponding to positive
terms early on in the Hilbert series, but this ceases to be logically
reliable once relations appear. For each singularity $\recip{r}(a,r-a)$
in the basket, there must be generators of weight divisible by $r$ to
cancel the periodicity in the Hilbert series, and a generator of
weight $\equiv a$ and $-a$ to provide local orbifold coordinates at
the singularity. We usually expect a generator in each degree $r$ as
in the above example, but sometimes a GCD of two different $r$ can
account for two at one go.

Deploying these tricks is something of an art, and there is no point
in being precise here about how we do it. Experience suggests that
they need to be combined in different ways in different situations (we
find different and better ways of putting in the generators
systematically each time we extend the database). But however we use
the numerical tricks on the Hilbert series, there will usually be work
left to do. Particular features of the geometry may demand extra
generators not predicted by the numerical data. For K3 surfaces in
higher codimension, the existence of particular linear systems often
forces us to include extra generators. One reason for stopping at
codimension~4 with the current version of the Magma database is that
there is little extra work in that case.

\subsubsection{{\tt MakeK3Database(BB,gg,cmax)}}

This function takes as its arguments a list of baskets $BB$ with list
of minimal genus $gg$, and an integer {\tt cmax}. It takes each basket
and genus pair in turn and applies any analysis of the Hilbert series
that is implemented. If the codimension gets larger than {\tt cmax},
the result is discarded (in the current implementation --- this is
handled in a more useful way in the proto\-type future implementation).
Once done, the result is wrapped in cosmetics and available for computer
study.

Of course, as we learn more tricks for analysing Hilbert series, we can
run them through the database, modifying candidate surfaces as we go.
But we believe that the current Magma database of 391 K3 surfaces is
reliable and complete. (For that matter, we are fairly confident of the
codimension~5 list in \cite{Br1}, which has $N_5=162$ elements.)

\subsubsection{{\tt Centres($\sim$DB)}, etc.}

There are several functions used to study the database by hand. The
procedure
\begin{verbatim}
> DB := K3Database("t");
> Centres(~DB); // Takes some minutes.
\end{verbatim}
performs the projection calculus of Exercise~\ref{ex!t1} to compute all
possible projections between the surfaces of the database {\tt DB}; the
tilde indicates that the command is {\em procedural\/}, that is, it
actually modifies {\tt DB} by writing in the centres it computes (the
next version of the database will have the centres ready for use on
loading). Once the centres are in, we can look at the projections from
any surface. For example, we can find the surface $S$ of \ref{exa!dbeg1}
in the database and then look for all possible iterated projections
from $S$:
\begin{verbatim}
> S := K3SurfaceFromWeights(DB,[3,4,5,6,7,10,13]);
> pc := ProjectionChains(S,DB); #pc;
4
> pc[1];
[
 [ [ 254, 4 ], [ 5, 1, 4 ], [ 10 ] ],
 [ [ 10, 1 ], [ 13, 3, 10 ], [ 40 ] ],
 [ [ 40, 1 ], [ 10, 3, 7 ], [ 79 ] ],
 [ [ 79, 1 ] ]
]
\end{verbatim}
This output is opaque (Leitmotif: please buy the upgrade, with its
many major enhancements), but you can probably see what is going on.
There are 4 chains of projections starting with $S$, and {\tt pc[1]}
only asks for the first. If we denote the $n$th surface in the database
by $S_n$, then $S=S_{254}$ and the first chain of projections, returned
by {\tt pc[1]}, is
\[
S\broken S_{10}\broken S_{40}\broken S_{79}.
\]
The output records the codimension of each surface. In this case the
final three surfaces are all codimension 1: the projections between
them are the quadratic involutions of \cite{CPR}. The second column
gives the type of centre {\tt [r,a,r-a]} of each projection.

The fourth chain of projections is
\begin{verbatim}
> pc[4];
[
 [ [ 254, 4 ], [ 13, 3, 10 ], [ 42 ] ],
 [ [ 42, 3 ], [ 10, 3, 7 ], [ 81 ] ],
 [ [ 81, 2 ], [ 7, 3, 4 ], [ 107 ] ],
 [ [ 107, 1 ] ]
]
\end{verbatim}
It consists purely of Type~I projections. In terms of homogeneous
coordinate rings, at each stage a single variable (of weight 13, 10, 7
respectively, the index of the corresponding centre) is eliminated from
the ring. A basic exercise in Magma makes this very clear:
\begin{verbatim}
> [ Weights(DB[i]) : i in [254,42,81,107] ];
[
    [ 3, 4, 5, 6, 7, 10, 13 ],
    [ 3, 4, 5, 6, 7, 10 ],
    [ 3, 4, 5, 6, 7 ],
    [ 3, 4, 5, 6 ]
]
\end{verbatim}
terminating with the surface $S_{107}$ which is
\begin{verbatim}
> DB[107];
Codimension 1 K3 surface, number 107, Reid1(39), with data
  Weights: [ 3, 4, 5, 6 ]
  Numerator: -t^18 + 1
  Basket: [ 2, 1 ], [ 3, 1 ], [ 3, 1 ], [ 3, 1 ], [ 4, 1 ], [ 5, 1 ]
\end{verbatim}
 \begin{exc} Write down a hypersurface $S_{18}\subset\PP(3,4,5,6)$
that contains $\PP(3,4)$ as the unprojection divisor of the final
projection, and singularities on it to provide the inverses of
the successive unprojections.
 \end{exc}


\begin{thebibliography}{CPR}

 \bibitem[A]{A} S. Alt{\i}nok, Graded rings corresponding to polarised
K3 surfaces and $\Q$-Fano 3-folds, Univ. of Warwick PhD thesis, Sep.
1998, $93+\mathrm{vii~pp.}$, get from
www.maths.warwick.ac.uk/\linebreak[2]$\!\sim$miles/doctors/Selma

 \bibitem[A1]{A1} S. Alt{\i}nok, Hilbert series and applications to
graded rings, submitted

 \bibitem[Be]{Be} Sarah-Marie Belcastro, Picard lattices of families of
K3 surfaces, to appear in Communications in Algebra, preprint
math.AG/9809008, 21 pp.

 \bibitem[Br]{Br} Gavin Brown, Datagraphs in algebraic geometry,
submitted to Proceedings of SNSC '01 (Linz), Ed. F. Winkler

 \bibitem[GRDW]{Br1} Gavin Brown, Magma algebraic geometry database
website, see
www.maths.\linebreak[2]warwick.\linebreak[2]ac.uk/$\!\sim$gavinb/grdb.html

 \bibitem[Bu]{Bu} Anita Buckley, Graded rings over polarised
Calabi--Yau 3-folds, work in progress

 \bibitem[CPR]{CPR} A. Corti, A. Pukhlikov and M. Reid, Birationally
rigid Fano hypersurfaces, in Explicit birational geometry of 3-folds,
A. Corti and M. Reid (eds.), CUP 2000, 175--258

 \bibitem[CR]{CR} A. Corti and M. Reid, Weighted Grassmannians, first
draft submitted to Proc. Francia memorial conference (Genova, 2001),
currently 12~pp.

 \bibitem[D]{D} I. Dolgachev, Weighted projective varieties, in Group
actions and vector fields (Vancouver, B.C., 1981), LNM{\bf956},
pp.~34--71

 \bibitem[Fl]{Fl} A.R. Iano-Fletcher, Working with weighted complete
intersections, in Explicit birational geometry of 3-folds, CUP 2000,
pp.~101--173

 \bibitem[GW]{GW} S. Goto and K-i Watanabe, On graded rings. I, J. Math.
Soc. Japan {\bf30} (1978) 179--213

 \bibitem[G]{G} M. Gross, Deforming Calabi-Yau threefolds, Math. Ann.
{\bf308} (1997) 187--220

 \bibitem[I1]{I1} V. A. Iskovskikh, Fano threefolds. I, Izv. Akad. Nauk
SSSR Ser. Mat. {\bf41} (1977) 516--562 = Math. USSR--Izv. {\bf11}
(1977) 485--527
 
 \bibitem[I2]{I2} V. A. Iskovskikh, Fano threefolds. II, Izv. Akad. Nauk
SSSR Ser. Mat. {\bf42} (1978) 506--549 = Math. USSR--Izv. {\bf12} (1978)
469--506

 \bibitem[Ka1]{Ka1} KAWAMATA Yujiro, On the plurigenera of minimal
algebraic 3-folds with $K\approx0$, Math. Ann. {\bf275} (1986) 539--546

 \bibitem[Ka2]{Ka2} KAWAMATA Yujiro, Boundedness of $\Q$-Fano
threefolds, in Proc. Internat. Conference on Algebra (Novosibirsk,
1989), Contemp. Math. {\bf131}, AMS, 1992, Part~3, pp.~439--445

 \bibitem[Ke]{Ke} Adam Keenan, Riemann--Roch, subcanonical divisors on
curves and graded rings, work in progress

 \bibitem[KM]{KM} A. Kustin and M. Miller, Constructing big Gorenstein
ideals from small ones, J. Algebra {\bf85} (1983) 303--322

 \bibitem[Leng]{Leng} R. Leng, McKay correspondence and equivariant
Riemann--Roch, Warwick PhD thesis in preparation

 \bibitem[Ma]{Ma} Magma (John Cannon's computer algebra system): W.
Bosma, J. Cannon and C.~Playoust, The Magma algebra system I: The user
language, J. Symb. Comp. {\bf24} (1997) 235--265. See also
www.maths.\linebreak[2]usyd.edu.au:8000/u/magma

 \bibitem[Mi1]{Mi1} MINAGAWA Tatsuhiro, Deformations of $\Q$-Calabi--Yau
3-folds and $\Q$-Fano 3-folds of Fano index~1, math.AG/9905106, J.
Math. Sci. Univ. Tokyo Tokyo {\bf6} (1999) 397-414

 \bibitem[Mi2]{Mi2} MINAGAWA Tatsuhiro, Deformations of weak Fano
3-folds with only terminal singularities, math.AG/9905107, Osaka J.
Math. {\bf38} (2001), to appear 

 \bibitem[MM1]{MM1} MORI Shigefumi and MUKAI Shigeru, Classification of
Fano 3-folds with $B_2\ge2$, Manuscripta Math. {\bf36} (1981/82)
147--162

 \bibitem[MM2]{MM2} MORI Shigefumi and MUKAI Shigeru, On Fano $3$-folds
with $B_2\ge2$, in Algebraic varieties and analytic varieties
(Tokyo, 1981), Adv. Stud. Pure Math. {\bf1}, Kinokuniya and
North-Holland, 1983, pp.~101--129

 \bibitem[Mu1]{Mu1} MUKAI Shigeru, Curves and Grassmannians, in
Algebraic geometry and related topics (Inchon, 1992), Internat. Press,
Cambridge MA 1993, pp.~19--40

 \bibitem[Mu2]{Mu2} Mukai Shigeru, Curves and symmetric spaces. I, 
Amer. J. Math. {\bf117} (1995) 1627--1644

 \bibitem[N]{N} NAMIKAWA Yoshinori, Smoothing Fano $3$-folds, J. alg
geom {\bf6} (1997) 307--324

 \bibitem[NS]{NS} NAMIKAWA Yoshinori and J.H.M. Steenbrink, Global
smoothing of Calabi-Yau threefolds, Invent. Math. {\bf122} (1995)
403--419
 
 \bibitem[P]{P} Stavros Papadakis, Gorenstein rings and Kustin--Miller
unprojection, Univ. of Warwick PhD thesis, Aug 2001, vi + 72~pp., get
from
www.maths.warwick.ac.uk/$\!\sim$miles/doctors/\linebreak[2]Stavros

 \bibitem[P1]{P1}
 Stavros Papadakis, Kustin--Miller unprojection {\em with} complexes,
submitted, math.AG/\linebreak[2]0111195, 23 pp.

 \bibitem[PR]{PR} Stavros Papadakis and Miles Reid, Kustin--Miller
unprojection without complexes, J. algebraic geometry (to appear),
preprint math.AG/0011094, 18~pp.

 \bibitem[YPG]{YPG} M. Reid, Young person's guide to canonical
singularities, in Algebraic geometry (Bowdoin, 1985), Proc. Sympos. Pure
Math. {\bf46} Part~1, AMS 1987, pp.~345--414,
 
 \bibitem[R1]{R1} M. Reid, Projective morphisms according to Kawamata,
Warwick preprint, 1983, unpublished, get from
www.maths.warwick.\linebreak[2]ac.uk/$\!\sim$miles/3folds

 \bibitem[Ki]{Ki} M. Reid, Graded rings and birational geometry, in
Proc. of algebraic geo\-metry symposium (Kinosaki, Oct 2000), K. Ohno
(Ed.), 1--72, get from www.maths.warwick.ac.uk/\linebreak[2]$\!\sim$miles/3folds

 \bibitem[S]{S} K. Suzuki, Univ. of Tokyo PhD thesis, in preparation

 \bibitem[T4]{T4} M. Reid, Examples of Type~IV unprojection, preprint
math.AG/\linebreak[2]0108037, 16~pp.

 \bibitem[T]{T} TAKAGI Hiromichi, On the classification of $\Q$-Fano
3-folds of Gorenstein index~2. I, II, RIMS preprint 1305, Nov 2000,
66~pp.

 \bibitem[T1]{T1} TAKAGI Hiromichi, work in progress

 \bibitem[W]{W} WATANABE Kei-ichi, Some remarks concerning Demazure's
construction of normal graded rings, Nagoya Math. J. {\bf 83} (1981)
203--211

 \end{thebibliography}
\end{document}